\documentclass[12pt,twoside]{amsart}
\usepackage{amssymb}

\usepackage[all]{xy}
\nonstopmode

\headsep=0.5cm  \numberwithin{equation}{section}
\font\tengothic=eufm10 scaled\magstep 1 \font\sevengothic=eufm7
scaled\magstep 1
\newfam\gothicfam
      \textfont\gothicfam=\tengothic
      \scriptfont\gothicfam=\sevengothic

\newtheorem{theorem}{Theorem}[section]
\newtheorem{lemma}[theorem]{Lemma}
\newtheorem{proposition}[theorem]{Proposition}
\newtheorem{corollary}[theorem]{Corollary}

\theoremstyle{definition}
\newtheorem{definition}[theorem]{Definition}
\newtheorem{remark}[theorem]{Remark}
\newtheorem{problem}[theorem]{Problem}
\newtheorem{example}[theorem]{Example}

\newtheorem{notation}[theorem]{Notation}

\newcommand{\Hom}{\operatorname{Hom}}

\newcommand{\Ext}{\operatorname{Ext}}

\newcommand{\cB}{{\mathcal B}}
\newcommand{\cE}{{\mathcal E}}
\newcommand{\cT}{{\mathcal T}}
\newcommand{\cF}{{\mathcal F}}
\newcommand{\cH}{{\mathcal H}}
\newcommand{\cG}{{\mathcal G}}
\newcommand{\cO}{{\mathcal O}}
\newcommand{\cD}{{\mathcal D}}
\newcommand{\cL}{{\mathcal L}}

\newcommand {\CC}{\mathbb{C}}

\newcommand {\ZZ}{\mathbb{Z}}

\newcommand {\PP}{\mathbb{P}}
\newcommand {\FF}{\mathbb{F}}
\newcommand {\mm}{{\textbf{m}}}
\newcommand {\mmm}{{\textrm{m}}}
\newcommand{\cS}{{\mathcal S}}
\newcommand{\cQ}{{\mathcal Q}}

\def\mapright#1{\smash{ \mathop{\longrightarrow}
    \limits^{#1}}}

\begin{document}
\title[\mmm-blocks collections and Castelnuovo-Mumford regularity ]
 {\mm-blocks collections and Castelnuovo-Mumford regularity in multiprojective spaces}

\author[L.\ Costa, R.M.\ Mir\'o-Roig]{L.\ Costa$^*$, R.M.\
Mir\'o-Roig$^{**}$}

\address{Facultat de Matem\`atiques,
Departament d'Algebra i Geometria, Gran Via de les Corts Catalanes
585, 08007 Barcelona, SPAIN } \email{costa@ub.edu}

\address{Facultat de Matem\`atiques,
Departament d'Algebra i Geometria, Gran Via de les Corts Catalanes
585, 08007 Barcelona, SPAIN } \email{miro@ub.edu}

\date{\today}
\thanks{$^*$ Partially supported by MTM2004-00666.}
\thanks{$^{**}$ Partially supported by MTM2004-00666.}



\begin{abstract}
 The main goal of the paper is to generalize Castelnuovo-Mumford
regularity for coherent sheaves on projective spaces to coherent
sheaves on $n$-dimensional smooth projective varieties $X$ with an
$n$-block collection $\cB $ which generates the bounded derived
category $\cD ^b({\cO}_X$-$mod)$. To this end, we use the theory
of $n$-blocks and Beilinson type spectral sequence to define the
notion of regularity of a coherent sheaf $F$ on $X$ with respect
to the $n$-block collection $\cB $. We show that the basic formal
properties of the Castelnuovo-Mumford regularity of coherent
sheaves over projective spaces continue to hold in this new
setting and we compare our definition of regularity with previous
ones. In particular, we show that in case of coherent sheaves on
$\PP^n$ and for the $n$-block collection $\cB
=(\cO_{\PP^n},\cO_{\PP^n} (1), \cdots , \cO_{\PP^n}(n))$ on
$\PP^n$ Castelnuovo-Mumford regularity and our new definition of
regularity coincide. Finally, we carefully study the regularity of
coherent sheaves on a multiprojective space $\PP^{n_1}\times
\cdots \times \PP^{n_r}$ with respect to a suitable $n_1+\cdots
+n_r$-block collection and we compare it with the multigraded
variant of the Castelnuovo-Mumford regularity given by Hoffman and
Wang in \cite{HW}.
\end{abstract}


\maketitle

\tableofcontents


 \section{Introduction} \label{intro}

In chapter 14 of \cite{M} D. Mumford introduced the concept of
regularity for a coherent sheaf $F$ on a projective space $\PP^n$
to bound the family of all projective subschemes having fixed
Hilbert polynomial. Since then Castelnuovo-Mumford regularity has
become a fundamental invariant in commutative algebra and
algebraic geometry. It measures the complexity of a module or a
sheaf; more precisely  the regularity of a module bounds the
largest degree of the minimal generators and the degree of
syzygies and the regularity of a sheaf estimates the smallest
twist for which the sheaf is globally generated.

\vskip 2mm Let $X$ be a smooth projective variety of dimension $n$
and let $\cB=(\cE_0, \cdots, \cE_n)$ be an $n$-block collection of
objects of  $\cD ^b({\cO}_X$-$mod)$  which generates the bounded
derived category  $\cD ^b({\cO}_X$-$mod)$. The goal of this paper
is to introduce the notion of regularity of a coherent sheaf on
$X$ with respect to $\cB$ as a generalization of the notion of
Castelnuovo-Mumford regularity of coherent sheaves on projective
spaces. To introduce this new notion of regularity and to state
the basic formal properties, we use helix theory, $m$-blocks
collections and Beilinson type spectral sequences.

\vskip 2mm We want to stress that Castelnuovo-Mumford regularity
as well as the notion of regularity developed for Grassmannians in
\cite{Ch} and multiprojective spaces in \cite{HW} fall under the
umbrella of $\cB $-regularity. Moreover,  in this new setting, we
are able to prove analogs of some of the classical results on
$m$-regularity for coherent sheaves on projective spaces.

\vskip 4mm Next we outline the structure of the paper. In section
2, we briefly recall the notions and properties of full strongly
exceptional collections of sheaves  on a smooth projective variety
needed later. The notion of $m$-block as well as the concept of
mutations of blocks  are presented in section 3. Sections 4 and 5
are the heart of the paper. In section 4, we  first introduce the
notion of helix of blocks associated to an $m$-block collection of
objects of $\cD ^b({\cO}_X$-$mod)$ as a natural generalization of
the notion of helix associated to an exceptional collection of
objects of $\cD ^b({\cO}_X$-$mod)$. Then, using Beilinson type
spectral  sequences, we give the promised definition of regularity
of a coherent sheaf $F$ on $X$ with respect to an $n$-block
collection $\cB $ which generates $\cD ^b({\cO}_X$-$mod)$, we
prove that the Castelnuovo-Mumford regularity of a coherent sheaf
$F$ on $\PP^n$ coincides with the regularity of $F$ with respect
to a suitable $n$-block collection on $\PP^n$ and we show that the
main formal properties of the Castelnuovo-Mumford regularity carry
over to the new setting. In section 5, we restrict our attention
to multiprojective spaces $X=\PP^{n_1}\times \cdots \times
\PP^{n_r}$ and we analyze the relationship between our definition
of regularity and the multigraded variant of the
Castelnuovo-Mumford regularity introduced by Hoffman and Wang in
\cite{HW}.
 Finally, in section 6, we
collect some questions which naturally arise from this paper.


\vskip 4mm  \noindent \underline{Notation}
 Throughout this paper $X$ will be a smooth projective variety defined over the
complex numbers $\CC$ (most of the results are true for varieties
over an algebraically closed field)  and we denote by $\cD=\cD
^b({\cO}_X$-$mod)$ the derived category of bounded complexes of
coherent sheaves of ${\cO}_X$-modules. Notice that $\cD$ is an
abelian linear triangulated category. We identify, as usual, any
coherent sheaf $F$ on $X$ to the object $(0 \rightarrow F
\rightarrow 0) \in \cD$ concentrated in degree zero and we will
not distinguish between a vector bundle and its locally free sheaf
of sections.


\section{Basic facts on exceptional collections}

As we pointed out in the introduction, in this section we gather
the basic definitions and properties on  full strongly exceptional
collections of sheaves  on a smooth projective variety needed in
the sequel. For general facts on triangulated categories see
\cite{V}.

\begin{definition}\label{exceptcoll}
Let $X$ be a smooth projective variety.

(i) An object  $F \in \cD$ is {\bf exceptional} if
$\Hom^{\bullet}_{\cD} (F,F)$ is a 1-dimensional algebra generated
by the identity.

(ii) An ordered collection $(F_0,F_1,\cdots ,F_m)$  of objects of
$\cD$ is an {\bf exceptional collection} if each object $F_{i}$ is
exceptional and $\Ext^{\bullet}_{\cD}(F_{k},F_{j})=0$ for $j<k$.

(iii) An exceptional collection $(F_0,F_1,\cdots ,F_m)$ of objects
of $\cD$ is a {\bf strongly exceptional collection} if in addition
$\Ext^{i}_{\cD}(F_j,F_k)=0$ for $i\neq 0$ and  $j \leq k$.

(iv) An ordered collection of objects of $\cD$, $(F_0,F_1,\cdots
,F_m)$,  is a {\bf full (strongly) exceptional collection} if it
is a (strongly) exceptional collection  and $F_0$, $F_1$, $\cdots
$ , $F_m$ generate the bounded derived category $\cD$.
\end{definition}

\begin{remark} The existence of a full strongly exceptional collection
$(F_0,F_1,\cdots,F_m)$ of coherent sheaves on a smooth projective
variety $X$ imposes rather a strong restriction on $X$, namely
that the Grothendieck group $K_0(X)=K_0(\cO _X-mod)$ is isomorphic
to $\ZZ^{m+1}$.
\end{remark}

Let us illustrate the above definition with precise examples:

\begin{example} \label{prihirse2} (1) ($\cO_{\PP^r}$, $\cO_{\PP^r} (1) $,
 $\cdots $, $\cO_{\PP^r}(r)$) is a full strongly exceptional
collection of coherent sheaves on  $\PP^r$ and ($\Omega^r_{\PP^r}
(r)$, $\Omega^{r-1}_{\PP^{r}} (r-1) $,  $\cdots
$,$\Omega^{1}_{\PP^{r}} (1) $,
 $\cO_{\PP^r}$) is also a full strongly
exceptional collection of coherent sheaves on $\PP^r$.

(2) Let $\FF_n=\PP(\cO_{\PP^1}\oplus \cO_{\PP^1}(n))$, $n\ge 0$,
be a Hirzebruch surface. Denote by $\xi$ (resp. $F$) the class of
the tautological line bundle (resp. the class of a fiber of the
natural projection $p: \FF_n \rightarrow \PP^1$). Then, ($\cO$,
$\cO (F) $, $\cO (\xi) $, $\cO (F+\xi) $) is a full strongly
exceptional collection of coherent sheaves on $\FF_n$.

(3) $( \cO_{\PP^n}(-n) \boxtimes \cO_{\PP^m}(-m),
\cO_{\PP^n}(-n+1) \boxtimes \cO_{\PP^m}(-m),
 \cdots, \cO_{\PP^n} \boxtimes \cO_{\PP^m}(-m),
\cdots,  \cO_{\PP^n}(-n) \boxtimes \cO_{\PP^m}, \cO_{\PP^n}(-n+1)
\boxtimes \cO_{\PP^m},$ $ \cdots, \cO_{\PP^n} \boxtimes
\cO_{\PP^m}) $ is a full strongly exceptional collection of
locally free sheaves on $\PP^n \times \PP^m$ (See also
\cite{CRMMR}; Proposition 4.16).

(4) Let $\pi: \widetilde{\PP}^2(l) \rightarrow \PP^2$  be the blow
up of $\PP^2$ at $l$ points and let $E_1=\pi^{-1}(p_1), \ldots,
E_l=\pi^{-1}(p_l)$ be the exceptional divisors. Then, the
collection of line bundles on $\widetilde{\PP}^2(l)$
\[ (\cO, \cO(E_1), \cO(E_2), \cdots, \cO(E_l), \cO(H),\cO(2H)) \]
 is a full strongly exceptional collection of coherent sheaves on $\widetilde{\PP}^2(l)$.

(5) Let  $X=Gr(k,n)$ be the Grassmannian of $k$-dimensional
subspaces of the $n$-dimensional vector space.  We have  the
canonical exact sequence
\[ 0 \rightarrow \cS \rightarrow \cO_X^n \rightarrow \cQ \rightarrow 0
\] where $\cS$ denotes the tautological $k$-dimensional bundle and
$\cQ$ the quotient bundle. In the sequel, $\Sigma^{\alpha}\cS$
denotes the space of the irreducible representations of the group
$GL(\cS)$ with highest weight $\alpha=(\alpha_1, \ldots,
\alpha_s)$ and $|\alpha|=\sum_{i=1}^{s} \alpha_i$ (see \cite{FH}
for general facts on Weyl modules). Denote by $A(k,n)$ the set of
locally free sheaves $\Sigma^{\alpha}\cS$ on $Gr(k,n)$ where
$\alpha$ runs over Young diagrams fitting inside a $k \times
(n-k)$ rectangle. Notice that for any $\Sigma^{\alpha} \cS \in
A(k,n)$, $0 \leq |\alpha| \leq k(n-k)$. Set $\rho(k,n):= \sharp
A(k,n)$. By \cite{Kap}; Proposition 2.2 (a), for any
$\Sigma^{\alpha} \cS, \Sigma^{\beta}\cS \in A(k,n)$,
$\Ext^i(\Sigma^{\alpha} \cS, \Sigma^{\beta}\cS) \neq 0$ only if
$i=0$ and by \cite{Kap2} (3.5), $\Hom(\Sigma^{\alpha} \cS,
\Sigma^{\beta}\cS) \neq 0$ only if $\alpha_i \geq \beta_i$ for all
$i$. Denote by $\cE _r$ the set of bundles $\Sigma^{\alpha} \cS
\in A(k,n)$ with $|\alpha|= k(n-k)-r$. Let $\sigma$ be the ordered
collection of locally free sheaves on $X$ constructed in the
following way. Going from the left to the right, first put all the
$\Sigma^{\alpha} \cS \in \cE_0$, i.e. all the $\Sigma^{\alpha} \cS
\in A(k,n)$ with $|\alpha|= k(n-k)$. The $i$-th time put all the
$\Sigma^{\alpha}\cS \in \cE_i$, i.e. all the $\Sigma^{\alpha} \cS
\in A(k,n)$ with $|\alpha|= k(n-k)-i$ and so on until $i=k(n-k)$.
By construction $\sigma$ is a strongly exceptional collection and
by \cite {Kap} Proposition 1.4, it is full. So, $A(k,n)$ can be
totally ordered in such a way that we obtain a full strongly
exceptional collection ($ E_1, \ldots, E_{\rho(k,n)}$) of locally
free sheaves on $X$.

(6)  Let  $Q_n \subset \PP^{n+1}$, $n>2$, be a hyperquadric
surface.  By \cite{Kap2}; Proposition 4.9, if $n$ is even and
$\Sigma _1$, $\Sigma _2$ are the Spinor bundles on $Q_n$, then
\[(\Sigma _1(-n), \Sigma _2(-n), \cO_{Q_n}(-n+1), \cdots, \cO_{Q_n}(-1), \cO_{Q_n}) \]
is a full strongly exceptional collection of locally free sheaves
on $Q_n$; and if $n$ is odd and $\Sigma$ is the Spinor bundle on
$Q_n$, then
\[( \Sigma(-n), \cO_{Q_n}(-n+1), \cdots, \cO_{Q_n}(-1), \cO_{Q_n}) \]
is a full strongly exceptional collection of locally free sheaves
on $Q_n$.
\end{example}

The importance of the existence of full strongly exceptional
collections relies on  the fact that  each full strongly
exceptional collection $(F_0,F_1,\cdots,F_m)$ of coherent sheaves
on a smooth projective variety $X$ determines a tilting sheaf
$\cT=\oplus_{i=0}^mF_{i}$ and hence functors
\textbf{R}$\Hom_X(\cT,-): D^b({\cO}_X-mod) \longrightarrow D^b(A)$
and $-\otimes_A^{\mbox{\textbf{L}}}\cT:D^b(A)\longrightarrow
D^b({\cO}_X-mod)$ which define mutually inverse equivalences
between the bounded derived categories of coherent sheaves on $X$
and the bounded derived category of finitely generated right
$A=\Hom_X(\cT,\cT)$-modules, respectively.

\begin{definition}
\label{mutation} Let $X$ be a smooth projective variety and let
$(A,B)$ be an exceptional pair of objects of $\cD$. We define
objects $L_AB$ and $R_BA$ with the aid of the following
distinguished triangles in the category $\cD$:
\begin{equation}  \label{t1} L_AB \rightarrow \Hom^{\bullet}_{\cD}(A,B) \otimes A \rightarrow B
\rightarrow L_AB[1]
\end{equation}
\begin{equation}  \label{t2} R_BA[-1] \rightarrow A \rightarrow \Hom^{\times \bullet}_{\cD}(A,B)
 \otimes B \rightarrow R_BA.
\end{equation}
A {\bf left mutation} of an exceptional pair $\sigma=(A,B)$ is the
pair
\[L_A \sigma=(L_AB,A)=(LB,A) \]
and a  {\bf right mutation} of an exceptional pair $\sigma=(A,B)$
is the pair
\[R_A \sigma=(B,R_BA)=(B,RA). \]
Lower indices will be omitted whenever this does not cause
confusion.
\end{definition}

\vspace{3mm}
\begin{definition}
Let $X$ be a smooth projective variety and let $\sigma=(E_0,
\cdots, E_m )$ be  an exceptional collection of objects of $\cD$.
A {\bf left mutation} (resp. {\bf right mutation}) of $\sigma$ is
defined as a mutation of a pair of adjacent objects in this
collection, i.e. for any $1 \leq i \leq m$ a left mutation $L_i$
replaces the $i$-th pair of consequent elements $(E_{i-1}, E_i)$
by its left mutation $(L_{E_{i-1}}E_i, E_{i-1} )$ and a right
mutation $R_i$ replaces the same pair of consequent elements
$(E_{i-1}, E_i)$ by its right mutation $(E_i, R_{E_{i}}E_{i-1})$:
\[ L_i \sigma=L_{E_{i-1}} \sigma= ( E_0, \cdots ,L_{E_{i-1}}E_i,
E_{i-1}, \cdots E_m )\]
\[ R_i \sigma=R_{E_{i-1}} \sigma= (E_0, \cdots ,E_i, R_{E_{i}}E_{i-1},
\cdots, E_m ).\]
\end{definition}

\vskip 2mm

\begin{notation} \label{composicio} Let $X$ be a smooth projective variety and let $\sigma=(F_0,
\cdots, F_m )$ be  an exceptional collection of objects of $\cD$.
It is convenient to agree that for any $0 \leq i,j \leq m$ and
$i+j \leq m$,
\[  R^{(j)}F_{i}=  R^{(j-1)}RF_{i}= R_{F_{i+j}} \cdots \cdots
  R_{F_{i+2}} R_{F_{i+1}} F_{i}=: R_{F_{i+j} \cdots \cdots
  F_{i+2} F_{i+1}} F_{i} \]

  \[ R^{(0)}_{F_{i-1}}\sigma= \sigma \quad \quad
 R^{(j)}_{F_{i-1}}\sigma=  R^{(j-1)}_{RF_{i-1}}(R_{F_{i-1}}\sigma) \]
 and similar notation for compositions of left mutations.
\end{notation}

\begin{remark}\label{patologia} (1) If $X$ is a smooth projective variety and  $\sigma=(F_0, \cdots,
F_m )$ is  an exceptional collection of objects of $\cD$, then any
mutation of $\sigma$ is an exceptional collection. Moreover, if
$\sigma$ generates the category $\cD$, then the mutated collection
also generates $\cD$.

 (2) In general, a mutation of a strongly exceptional collection is
not a  strongly exceptional collection. In fact, take $X= \PP^1
\times \PP^1$ and consider the full strongly exceptional
collection $\sigma=(\cO_X, \cO_X(1,0), \cO_X(0,1), \cO_X(1,1) ) $
of line bundles on $X$. It is not difficult to check that the
mutated collection $$ (\cO_X, \cO_X(1,0),
L_{\cO_X(0,1)}\cO_X(1,1), \cO_X(0,1) )=(\cO_X, \cO_X(1,0),
\cO_X(-1,1), \cO_X(0,1) )$$ is no more a strongly exceptional
collection of line bundles on $X$.
\end{remark}

Let $X$ be a smooth projective variety of dimension $n$. It is
well known that if full strongly exceptional collections of
coherent sheaves on $X$ exist then all of them have the same
length and it is equal to the rank of $K_0(X)$. Even more, this
length is bounded below by $n+1$ because for any smooth projective
variety $X$ of dimension $n$  we have $rank(K_0(X)) \geq n+1$. In
\cite{CMR}; we give the following definition (see also \cite{BP}
and \cite{H}):

\vspace{3mm}

\begin{definition}
\label{defexcellent}  Let $X$ be a smooth projective variety of
dimension $n$. We say that an ordered collection of coherent
sheaves $\sigma=(E_0, \cdots, E_n)$ is a {\bf geometric
collection} if it is a full exceptional collection of coherent
sheaves on $X$ of minimal length, $n+1$, i.e. of length one
greater than the dimension of $X$.
\end{definition}

\vspace{3mm}

By  \cite{Bo}; Assertion 9.2, Theorem 9.3 and Corollary 9.4,
geometric collections are automatically strongly exceptional
collections of coherent sheaves and  the strongly exceptionality
is preserved under mutations.

\begin{example}
\label{exemplesexcellents} (1)  The collection
$\sigma=(\cO_{\PP^r}(-r), \cO_{\PP^r} (-r+1) , \cO_{\PP^r} (-r+2)
, \cdots , \cO_{\PP^r} )$ of line bundles on $\PP^r$ is a
geometric collection of coherent sheaves.

(2) If $n$ is odd and $Q_n \subset \PP^{n+1}$ is a quadric
hypersurface,  the collection of locally free sheaves
\[( \Sigma(-n), \cO_{Q_n}(-n+1), \cdots, \cO_{Q_n}(-1), \cO_{Q_n}) \]
being $\Sigma$ the Spinor bundle on $Q_n$ is a geometric
collection of locally free sheaves on $Q_n$.

(3) If $n$ is even and $Q_n \subset \PP^{n+1}$ is a quadric
hypersurface,  the collection of locally free sheaves
\[(\Sigma _1(-n), \Sigma _2(-n), \cO_{Q_n}(-n+1), \cdots, \cO_{Q_n}(-1), \cO_{Q_n}) \]
being $\Sigma _1$ and $\Sigma _2$ the Spinor bundles on $Q_n$, is
a full strongly exceptional collection of locally free sheaves on
$Q_n$. Since all full strongly exceptional collections of coherent
sheaves on $Q_n$ have  length  $n+2$, we conclude that there are
no geometric collections of coherent sheaves on $Q_n$ for even
$n$.

(4) It follows from Example \ref{prihirse2} (5) that there are no
geometric collections of coherent sheaves on $Gr(k,n)$ if $k \neq
n-1$.

(5) Any smooth Fano threefold $X$ with $Pic(X) \cong \ZZ$ and
trivial intermediate Jacobian has a geometric collection (see
\cite{CMR}; Proposition 3.6).
\end{example}

\vskip 2mm In \cite{CMR}, the authors extend the notion of
Castelnuovo-Mumford regularity for coherent sheaves on projective
spaces to coherent sheaves on smooth projective varieties with a
geometric collection. So we are led to consider the following
problem:

\vspace{3mm}

\begin{problem}
 \label{prob1} To characterize the smooth
projective varieties which have a geometric collection.
\end{problem}

\vspace{3mm}

To our knowledge Problem \ref{prob1} is far of being solved (See
\cite{CMR} for more information). Moreover, we want to stress that
the existence of a geometric collection on an $n$-dimensional
smooth variety $X$ imposes a strong restriction on $X$; e.g. $X$
has to be a Fano variety (\cite{BP}; Theorem 3.4) and the
Grothendieck group $K_0(X)$ has to be a $\ZZ $-free module of rank
$n+1$. So, it is convenient to generalize the notion of geometric
collection in order to be able to extend the concept of
Castelnuovo-Mumford regularity for coherent sheaves on projective
spaces to coherent sheaves on smooth projective varieties as
Grassmannians, even-dimensional hyperquadrics, multiprojective
spaces, etc, which do not have geometric collections. This will be
achieved allowing exceptional collections $\sigma=(F_0,
\ldots,F_m)$ of arbitrary length but packing the objects $F_i \in
\cD$ in suitable subcollections called blocks.


\section{$m$-Blocks and mutations}

\vskip 2mm The notion of block was introduced by Karpov and Nogin
in \cite{KN}. We start this section recalling its definition and
properties (see also \cite{H}).

\vspace{3mm}

\begin{definition}
\label{block} (i) An exceptional collection $(F_0,F_1, \cdots,
F_m)$ of objects of $\cD$ is a {\bf block} if
$\Ext^{i}_{\cD}(F_j,F_k)=0$ for any $i$ and  $j \neq k$.

(ii) An {\bf $m$-block collection of  type} $(\alpha_0, \alpha_1,
\cdots, \alpha_m)$ of objects of $\cD$ is an exceptional
collection
\[\cB=(\cE_0, \cE_1, \cdots, \cE_m)=(E_1^0,\cdots,E_{\alpha_0}^0, E_1^1,
 \cdots, E_{\alpha_1}^1,
\cdots, E_1^m,\cdots, E_{\alpha_m}^m) \] such that all the
subcollections $\cE_j=(E_1^j,E_2^j,\cdots, E_{\alpha_j}^j)$ are
blocks.
\end{definition}

Note that an exceptional collection $(E_0,E_1, \cdots, E_m)$ is an
$m$-block of type $(1,1, \cdots, 1)$.

\vspace{3mm}

\begin{example}
\label{exempleblock} (1) ($\cO_{\PP^r}(-r)$, $\cO_{\PP^r} (-r+1)
$,
 $\cO_{\PP^r} (-r+2)
$, $\cdots $, $\cO_{\PP^r}$) is  an $r$-block of type $(1,1,
\cdots, 1)$.

 \vskip 2mm (2)  Let  $X=Gr(k,n)$ be the Grassmannian of
$k$-dimensional subspaces of the $n$-dimensional vector space,
$k>1$. In Example \ref{prihirse2} (4), we have seen that $A(k,n)$
can be totally ordered in such a way that we obtain a full
strongly exceptional collection \[ \sigma=( E_1, \ldots,
E_{\rho(k,n)})\] of locally free sheaves on $X$. Moreover, packing
in the same block $\cE_r$ the bundles $\Sigma^{\alpha} \cS \in
\sigma$ with $|\alpha|= k(n-k)-r$ we obtain
\[ \sigma=( E_1, \ldots, E_{\rho(k,n)})=(\cE_0, \ldots, \cE_{k(n-k)}) \]
a $k(n-k)$-block collection of locally free sheaves on $X$ (see
Example \ref{prihirse2} (4) for details).

\vskip 2mm (3) Let $Q_n \subset \PP^{n+1}$, $n \geq 2$, be a
hyperquadric variety. According to Example \ref{prihirse2} (5), if
$n$ is even and $\Sigma _1$, $\Sigma _2$ are the Spinor bundles on
$Q_n$, then
\[(\Sigma _1(-n), \Sigma _2(-n), \cO_{Q_n}(-n+1), \cdots, \cO_{Q_n}(-1), \cO_{Q_n}) \]
is a full strongly exceptional collection of locally free sheaves
on $Q_n$; and if $n$ is odd and $\Sigma$ is the Spinor bundle on
$Q_n$, then
\[( \Sigma(-n), \cO_{Q_n}(-n+1), \cdots, \cO_{Q_n}(-1), \cO_{Q_n}) \]
is a full strongly exceptional collection of locally free sheaves
on $Q_n$. Since $\Ext^i(\Sigma_1,\Sigma_2)=0$ for any $i \geq 0$,
we get that $(\cE_0, \cE_1, \ldots, \cE_n)$ where  \[
\cE_j=\cO_{Q_n}(-n+j) \quad \mbox{for } 1 \leq j \leq n, \quad
\cE_0=
\begin{cases} (\Sigma_1(-n), \Sigma_2(-n))  & \mbox{if } n
\quad \mbox{even}
\\ (\Sigma(-n)) & \mbox{if } n \quad \mbox{odd} \\
\end{cases} \] is an $n$-block collection of locally free sheaves on $Q_n$ for all
$n$.

\vskip 2mm (4) Let $X=\PP^{n_1} \times \cdots \times \PP^{n_s}$ be
a multiprojective space of dimension $d=n_1+\cdots+n_s$. For any
$1 \leq i \leq s$, denote by $p_i:X \rightarrow \PP^{n_i}$ the
natural projection and write
 \[ \cO_X(a_1, a_2, \cdots, a_s):= p_1^* \cO_{\PP^{n_1}}(a_1)
 \otimes p_2^*\cO_{\PP^{n_2}}(a_2) \otimes \cdots \otimes p_s^* \cO_{\PP^{n_s}}(a_s).\]
 For any $0 \leq j \leq d$, denote by $\cE_j$
 the collection of all line bundles on $X$
 \[ \cO_X(a_1^j, a_2^j, \cdots, a_s^j)\]
 with $-n_i \leq a_i^j \leq 0$ and
 $\sum_{i=1}^{s}a_i^j=j-d$.  Using the K\"{u}nneth formula for
 locally free sheaves on algebraic varieties, we prove that
 each $\cE_j$ is a block and that
 \[\cB=(\cE_0, \cE_1, \cdots, \cE_{d}) \]
 is a $d$-block collection of line bundles
 on $X$.

\end{example}

\vspace{3mm}

We will now introduce the notion of mutation of block collections.

\vspace{3mm}

\begin{definition} \label{mutationblock}
Let $X$ be a smooth projective variety and consider
  a 1-block
collection $(\cE,\cF)=(E_1, \cdots, E_n,F_1, \cdots, F_m)$ of
objects of $\cD$. A {\bf left mutation} of $F_j$ by $\cE$ is the
object defined by (see Notation \ref{composicio})
\[L_{\cE}F_j:=L_{E_1E_2 \cdots E_n}F_j \]
and a  {\bf right mutation} of $E_j$ by $\cF$ is the object
defined by
\[R_{\cF}E_j:= R_{F_mF_{m-1} \cdots F_1}E_j. \]
A {\bf left mutation} of $(\cE, \cF)$ is the pair $(L_{\cE}\cF,
\cE)$ where
\[L_{\cE}\cF:=(L_{\cE}F_1, L_{\cE}F_2, \cdots,L_{\cE}F_m) \]
and a {\bf right mutation} of $(\cE, \cF)$ is the pair
$(\cF,R_{\cF}\cE)$ where
\[R_{\cF}\cE:=(R_{\cF}E_1, R_{\cF}E_2, \cdots, R_{\cF}E_n). \]
\end{definition}

Note that by \cite{GK} (2.2),  $L_{\cE}\cF$ and $R_{\cF}\cE$ are
blocks and the pairs $(L_{\cE}\cF,\cE)$ and $(\cF, R_{\cF}\cE)$
are 1-block collections.

\begin{remark} \label{triangles}
It follows from the proof of \cite{KN}; Proposition 2.2 and
Proposition 2.3 that given a 1-block collection $(\cE,\cF)=(E_1,
\cdots, E_n,F_1, \cdots, F_m)$, the objects $L_{\cE}F_j$ and
$R_{\cF}E_j$ can be defined with the aid of the following
distinguished triangles in the category $\cD$
\begin{equation}
\label{triangleleft} L_{\cE}F_j \rightarrow \oplus_{i=1}^{n}
\Hom^{\bullet}_{\cD}(E_i,F_j) \otimes E_i \rightarrow F_j
\rightarrow L_{\cE}F_j[1]
\end{equation}
\begin{equation}
\label{triangleright} R_{\cF}E_j[-1] \rightarrow E_j \rightarrow
\oplus_{i=1}^{m} \Hom^{\times \bullet}_{\cD}(E_j,F_i)\otimes F_i
\rightarrow R_{\cF}E_j.
\end{equation}
\end{remark}

\vspace{3mm}

Applying $\Hom^{\bullet}_{\cD}(E_i, *)$ to the triangle
(\ref{triangleleft}) we get the orthogonality relation

\begin{equation}
\label{orto1l} \Hom^{\bullet}_{\cD}(E_i, L_{\cE} F_j)=0 \quad
\mbox{for all} \quad 1 \leq i \leq n
\end{equation}

\vskip 2mm \noindent i.e., $L_{\cE}F_j \in [\cE]^{\bot}:= \{F \in
\cD | \Hom^{\bullet}_{\cD}(E,F)=0 \quad \mbox{for all } E \in
[\cE]\}$, where we denote by $[\cE]$ the full triangulated
subcategory of $\cD$ generated by $E_1, \cdots, E_n$.

Similarly, $\Hom^{\bullet}_{\cD}(*, F_j)$ applied to the triangle
(\ref{triangleright}) gives the orthogonality relation

\begin{equation}
\label{orto1r} \Hom^{\bullet}_{\cD}( R_{\cF} E_i, F_j)=0 \quad
\mbox{for all} \quad 1 \leq j \leq m
\end{equation}

\vskip 2mm \noindent i.e., $R_{\cF}E_i \in$  $^{\bot}[\cF]:= \{E
\in \cD | \Hom^{\bullet}_{\cD}(E,F)=0 \quad \mbox{for all } F \in
[\cF]\}$.

\vspace{3mm}

\begin{notation} It is convenient to agree that
 \[  R^{(j)}\cE_{i}=  R^{(j-1)}R\cE_{i}= R_{\cE_{i+j}} \cdots \cdots
  R_{\cE_{i+2}} R_{\cE_{i+1}} \cE_{i}=: R_{\cE_{i+j} \cdots \cdots
  \cE_{i+2} \cE_{i+1}} \cE_{i} \]
   \[  L^{(j)}\cE_{i}=  L^{(j-1)}L\cE_{i}= L_{\cE_{i-j}} \cdots \cdots
  L_{\cE_{i-2}} L_{\cE_{i-1}} \cE_{i}=: L_{\cE_{i-j} \cdots \cdots
  \cE_{i-2} \cE_{i-1}} \cE_{i}. \]
\end{notation}

Let $\cB=(\cE_0, \cdots, \cE_m)$ be an $m$-block collection of
type $\alpha_0, \cdots, \alpha_m$ of objects of $\cD$ which
generates $\cD$. Two $m$-block collections $\cH=(\cH_0, \cdots,
\cH_m)$ and $\cG=(\cG_0, \cdots, \cG_m)$ of type  $\beta_0,
\cdots, \beta_m$ with $\beta_i=\alpha_{m-i}$ of objects of $\cD$
are called {\bf left dual $m$-block collection of $\cB $} and {\bf
right dual $m$-block collection of $\cB$} if
\begin{equation}\label{orto1}  \Hom^{\bullet}_{\cD}(H_{j}^{i},E^k_l)=\Hom^{\bullet}_{\cD}(E^k_l,
G^i_j)=0 \end{equation} except for

\begin{equation}\label{orto2}
 \Ext^{k}_{\cD}(H_{i}^{k},E_i^{m-k})=\Ext^{m-k}_{\cD}(E_i^{m-k},
G_i^k)=\CC.
\end{equation}

\begin{remark}\label{existence-dual}
 Let $X$ be a smooth projective variety. Given an $m$-block
collection $\cB=(\cE_0, \cdots, \cE_m)$ of type $\alpha_0, \cdots,
\alpha_m$ of objects of $\cD$ which generates $\cD$, left dual
$m$-block collections and right dual $m$-block collections  of
$\cB$ exist and they are unique up to isomorphism. In fact, by
\cite{CMR2}; Proposition 3.9, the $m$-block collection
\begin{equation} \label{caracdelsr} \cH=(R^{(0)}\cE_m,
R^{(1)}\cE_{m-1}, \cdots, R^{(m)}\cE_0)
\end{equation} where by definition
\[ \begin{array}{rl} R^{(i)}\cE_{m-i}= & (R^{(i)}E_1^{m-i}, \cdots,
R^{(i)}E_{\alpha_{m-i}}^{m-i} ) \\ = & (R_{\cE_m\cE_{m-1} \cdots
\cE_{m-i+1} }E_1^{m-i}, \cdots, R_{\cE_m\cE_{m-1} \cdots
\cE_{m-i+1} }E_{\alpha_{m-i}}^{m-i} ) \end{array} \]
 satisfies the orthogonality conditions $(\ref{orto1})$ and
$(\ref{orto2})$. Therefore, $\cH$ is  the left dual $m$-block
collection of $\cB $. By consequent left mutations of the
$m$-block collection $\cB $ and arguing in the same way we get the
right dual $m$-block collection of $\cB $.
\end{remark}

\vskip 3mm Let $X$ be an $n$-dimensional smooth projective variety
with an $m$-block collection $\cB=(\cE_0, \cdots, \cE_m)$ which
generates the bounded derived category $\cD$. The left dual
$m$-block collection of $\cB$ will play an important role in our
definition of regularity of a coherent sheaf $F$ on $X$ with
respect to $\cB$ (See Definition \ref{new}). Therefore, we will
now describe explicitly the left dual $m$-block collection of the
examples of $m$-block collections given in Example
\ref{exempleblock}.

\begin{example}\label{exampledualbase}
(1) Let $V$ be a $\CC $-vector space of dimension $n+1$ and set
$\PP^n=\PP(V)$. We consider the $n$-block collection  $\cB
=(\cO_{\PP^n}, \cO_{\PP^n}(1),\cdots ,\cO_{\PP^n}(n))$ on $\PP^n$.
Using the exterior powers $$ 0\longrightarrow \wedge
^{k-1}T_{\PP^n}\longrightarrow \wedge ^{k}V \otimes
\cO_{\PP^n}(k)\longrightarrow \wedge ^{k}T_{\PP^n}\longrightarrow
0$$ of the Euler sequence $$ 0\longrightarrow
\cO_{\PP^n}\longrightarrow V \otimes \cO_{\PP^n}(1)\longrightarrow
T_{\PP^n}\longrightarrow 0$$ we compute the left dual $n$-block
collection of $\cB =(\cO_{\PP^n}, \cO_{\PP^n}(1),\cdots
,\cO_{\PP^n}(n))$ and we get
$$(\cO_{\PP^n}(n), R^{(1)}\cO_{\PP^n}(n-1),\cdots
,R^{(j)}\cO_{\PP^n}(n-j),\cdots ,R^{(n)}\cO_{\PP^n})$$
$$=(\cO_{\PP^n}(n), T_{\PP^n}(n-1),\cdots
,\wedge^{j}T_{\PP^n}(n-j),\cdots ,\wedge^{n}T_{\PP^n}).$$

\vskip 2mm (2) Let $X=\PP^{n_1} \times \cdots \times \PP^{n_s}$ be
a multiprojective space of dimension $d=n_1+\cdots+n_s$.
 For any $0 \leq j \leq d$, denote by $\cE_j$
 the collection of all line bundles on $X$
 \[ \cO_X(a_1^j, a_2^j, \cdots, a_s^j)\]
 with $-n_i \leq a_i^j \leq 0$ and
 $\sum_{i=1}^{s}a_i^j=j-d$.  By Example \ref{exempleblock} (4),
 $\cB=(\cE_0, \cE_1, \cdots, \cE_{d}) $
 is a $d$-block collection of line bundles
 on $X$. By Remark \ref{existence-dual}, the left
dual $d$-block collection of $\cB$ is given by
\begin{equation} \label{caracdelsr} \cH=(R^{(0)}\cE_d,
R^{(1)}\cE_{d-1}, \cdots, R^{(d)}\cE_0)
\end{equation} where by definition
\[ \begin{array}{rl} R^{(k)}\cE_{d-k}= & (
\cdots,R^{(k)}\cO_{X}(t_1,\cdots,t_s), \cdots
 ) \\ = & ( \cdots ,  R_{\cE_d \cdots
 \cE_{d-k+1}}\cO_{X}(t_1,\cdots,t_s),
 \cdots
 ). \end{array} \]
 A straightforward computation shows that
  for any $\cO_{X}(t_1,\cdots,t_s) \in \cE_{d-k}$ and  any $0 \leq k \leq d$,
\[R^{(k)}\cO_{X}(t_1,\cdots,t_s)=
 R_{\cE_d \cdots \cE_{d-k+1}}\cO_{X}(t_1,\cdots,t_s)=
  \bigwedge^{-t_1} T_{\PP^{n_1}}(t_1)
\boxtimes \cdots \boxtimes \bigwedge^{-t_s} T_{\PP^{n_s}}(t_s).\]

 \vskip 2mm (3)  Let $Q_n \subset \PP^{n+1}$, $n \geq 2$, be a
hyperquadric variety and let  $\cB=(\cE_0, \cE_1, \ldots, \cE_n)$
where  \[ \cE_j=\cO_{Q_n}(-n+j) \quad \mbox{for } 1 \leq j \leq n,
\quad \cE_0=
\begin{cases} (\Sigma_1(-n), \Sigma_2(-n))  & \mbox{if } n
\quad \mbox{even}
\\ (\Sigma(-n)) & \mbox{if } n \quad \mbox{odd} \\
\end{cases} \] be the $n$-block collection of locally free sheaves on $Q_n$
 described in Example \ref{exempleblock}(3). To define the left
 dual $n$-block collection of $\cB$ we need to fix some notation.
 We set $\Omega ^{j}:= \Omega
^j_{\PP^{n+1}}$ and we define inductively $\psi _j$: $$\psi
_0:=\cO_ {Q_n}, \quad \psi _1:=\Omega^1(1)_{|Q_n}$$ and, for all
$j\ge 2$, we define the locally free sheaf $\psi _j$ as the unique
non-splitting extension (Note that $\Ext^1(\psi _{j-2},\Omega
^{j}(j)_{|Q_n})=\CC$): $$0\longrightarrow \Omega
^{j}(j)_{|Q_n}\longrightarrow \psi _j\longrightarrow \psi _{j-2}
\longrightarrow 0. $$ By \cite{Kap2}; Proposition 4.11 and using
the fact that the left dual $n$-block collection of a given
$n$-block collection  is uniquely determined up to unique
isomorphism by the orthogonality conditions described in Remark
\ref{existence-dual}, we get that the  left dual $n$-block
collection  of the $n$-block collection $\cB$ is
$$ \cH=(R^{(0)}\cE_n, R^{(1)}\cE_{n-1}, \cdots, R^{(n)}\cE_0)$$
where

 \[ R^{(j)}\cE_{n-j}=\psi_{j}^{*} \quad \mbox{for } 1 \leq j \leq n,
\quad R^{(n)}\cE_0=
\begin{cases} (\Sigma_1^{*}(1), \Sigma_2^{*}(1))  & \mbox{if } n
\quad \mbox{even}
\\ (\Sigma ^*(1)) & \mbox{if } n \quad \mbox{odd}. \\
\end{cases} \]

 \vskip 2mm (4) Let  $X=Gr(k,n)$ be the Grassmannian of
$k$-dimensional subspaces of the $n$-dimensional vector space and
let $\cB=(\cE_0,\cE_1, \cdots ,\cE_{k(n-k)})$ with $\cE _r=\{
\Sigma^{\alpha} \cS \mid |\alpha|= k(n-k)-r \}$ be the
$k(n-k)$-block collection of locally free sheaves on $X$ described
in Example \ref{exempleblock} (2). The sequence $\alpha $ defines
a Young diagram and we denote by $\tilde{\alpha }$ the sequence
corresponding to the conjugate diagram. By \cite{Kap2} Lemma 3.2,
for any two indices: $\alpha: $ $n-k\ge \alpha_1\ge \cdots \ge
\alpha_k\ge 0$ and $\beta:$ $k\ge \beta_1\ge \cdots \ge
\beta_{n-k}\ge 0$ we have $$H^{i}(X,\Sigma ^{\alpha}\cS\otimes
\Sigma ^{\beta}\cQ^*)=\begin{cases} \CC  & \mbox{if } \alpha
=\tilde{\beta} \quad \mbox{ and } i=| \alpha |
\\ 0 & \mbox{ otherwise. } \\
\end{cases}$$

Therefore, the left dual $k(n-k)$-block collection of
$\cB=(\cE_0,\cE_1, \cdots ,\cE_{k(n-k)})$ with $\cE _r=\{
\Sigma^{\alpha} \cS \mid |\alpha|= k(n-k)-r \}$ is

$$ \cH=(R^{(0)}\cE_{k(n-k)}, R^{(1)}\cE_{k(n-k)-1}, \cdots, R^{(j)}\cE_{k(n-k)-j}, \cdots,
 R^{(k(n-k))}\cE_0)
$$ where $ R^{(r)}\cE_{k(n-k)-r}=\{\Sigma ^{\tilde{\alpha}}\cQ
\mid |\alpha |=r \}.$
\end{example}

\vspace{3mm}

We want to point out that the notion of $m$-block collection is
the convenient generalization of the notion of geometric
collection we were looking for. Indeed, we will see that the
behavior of $n$-block collections, $n=\dim(X)$, is really good in
the sense that they are automatically strongly exceptional
collections and that their structure is preserved under mutations
through blocks. More precisely we have:

\begin{proposition}
\label{bonespropietats} Let $X$ be a smooth projective variety of
dimension $n$ and let $\cB=(\cE_0, \cdots, \cE_n)$ be an $n$-block
collection of coherent sheaves on $X$ and assume that $\cB$
generates the category $\cD$. Then, we have:

(1) The sequence $\cB $ is a full strongly exceptional collection
of coherent sheaves on $X$.

(2) All mutations through the blocks $\cE_j$ can be computed using
short exact sequences of coherent sheaves.

(3) Any mutation of $\cB$ through any block $\cE_j$ is a full
strongly exceptional collection of pure sheaves, i.e. complexes
concentrated in the zero component of the grading.

(4) Any mutation of $\cB$ through any block $\cE_j$ is an
$n$-block collection.
\end{proposition}
\begin{proof} See \cite{Bo}; Theorem 9.5 and Remark b) below and
\cite{H}; Theorem 1.
\end{proof}

\vspace{4mm}

\begin{remark} By Remark  \ref{patologia} (2), a mutation of a strongly exceptional collection is
not, in general,  a strongly exceptional collection. In fact,
$\sigma=(\cO_X, \cO_X(1,0), \cO_X(0,1), \cO_X(1,1) ) $ is a
strongly exceptional collection on $X=\PP^1 \times \PP^1$ and the
mutated collection $$ (\cO_X, \cO_X(1,0),
L_{\cO_X(0,1)}\cO_X(1,1), \cO_X(0,1) )=(\cO_X, \cO_X(1,0),
\cO_X(-1,1), \cO_X(0,1) )$$ is no more a strongly exceptional
collection on $X$. However, we can pack the objects of $\sigma $
in a suitable subcollections of blocks $\cB=(\cE_0,\cE_1,\cE_2)=
(\cO_X, (\cO_X(1,0), \cO_X(0,1)), \cO_X(1,1) )$ and according to
Proposition \ref{bonespropietats} any mutation of $\cB$ through
any block $\cE_j$ is a full strongly exceptional collection. So,
for instance, the mutation of $\cB$ through the block $\cE_1$ is
the full strongly exceptional collection
$$(\cE_0,L_{\cE_1}\cE_2,\cE_1)=(\cO_X,
L_{\cO_X(1,0)\cO_X(0,1)}\cO_X(1,1),(\cO_X(1,0), \cO_X(0,1))=$$ $$
(\cO_X, T_X(-1,-1),(\cO_X(1,0), \cO_X(0,1))).$$

To compute $L_{\cE_1}\cE_2=
L_{\cO_X(1,0)\cO_X(0,1)}\cO_X(1,1)=T_X(-1,-1)$ we have used the
exact sequences

$$ 0\longrightarrow \cO_X(-1,1) \longrightarrow V^{*}\otimes \cO_X
(0,1) \longrightarrow \cO_X(1,1) \longrightarrow 0, \mbox{ and }$$
$$ 0\longrightarrow \cO_X(-1,1) \longrightarrow T_X(-1,-1)
\longrightarrow \wedge ^2V^{*}\otimes \cO_X(1,0) \longrightarrow
0$$

\noindent being $X=\PP^1 \times \PP^1$ and $\PP^1=\PP(V)$.
\end{remark}

 \vspace{6mm}

\vspace{4mm} Beilinson Theorem was stated in 1978 \cite{Be} and
since then it has became a major tool in classifying vector
bundles over projective spaces. Beilinson spectral sequence was
generalized by Kapranov to hyperquadrics and Grassmannians
(\cite{Kap} and \cite{Kap2}) and by the authors to any smooth
projective variety with a geometric collection \cite{CMR}. We are
now ready to generalize Beilinson Theorem to any smooth projective
variety $X$ of dimension $n$
 with an $n$-block collection $\cB =(\cE_0, \cE_1,\cdots ,\cE _{n})$,
  $\cE_j=(E_1^j, \ldots,E_{\alpha_j}^j )$ of coherent sheaves on $X$
  which generates $\cD$.

 \vspace{4mm}

\begin{theorem} \label{mbe} {\bf (Beilinson type spectral sequence)}
Let $X$ be a smooth projective  variety of dimension $n$
 with an $n$-block collection $\cB =(\cE_0, \cE_1,\cdots ,\cE _{n})$,
  $\cE_i=(E_1^i, \ldots,E_{\alpha_i}^i )$ of coherent sheaves on $X$ which generates $\cD$.
   Then for any coherent sheaf $F$ on $X$ there are two spectral sequences
   situated in the square $-n \le p  \le 0$, $0 \le q \le n $, with $E_1$-term
 \begin{equation}\label{sucespectral0}
   _{I} E_1^{pq}=
\begin{cases}
   \bigoplus_{i=1}^{\alpha_{p+n}}\Ext^q(R_{\cE_{n} \cdots \cE_{p+n+1}}E _{i}^{p+n},{F})
   \otimes E_{i}^{p+n}   & \mbox{if} \quad -n \leq p \leq -1 \\
\bigoplus_{i=1}^{\alpha_{n}}\Ext^q(E _{i}^{n},{F})
   \otimes E_{i}^{n}   & \mbox{if} \quad p=0 \end{cases}
   \end{equation}
 \begin{equation}\label{sucespectral1}
   _{II} E_1^{pq}= \begin{cases} \bigoplus_{i=1}^{\alpha_{p+n}}\Ext^q( (E_{i}^{p+n})^*,{F})
   \otimes (R_{\cE_{n} \cdots \cE_{p+n+1}}E _{i}^{p+n})^* &
   \mbox{if} \quad -n \leq p \leq -1 \\
\bigoplus_{i=1}^{\alpha_{n}}\Ext^q( {E_{i}^{n}}^*,{F})
   \otimes {E _{i}^{n}}^* &
   \mbox{if} \quad  p=0 \end{cases}  \end{equation}
   and differentials $d_r^{pq}: E_r^{p,q} \rightarrow E_r^{p+r,q-r+1}
   $   which converge to $$ _{I} E^{i}_{\infty }=_{II} E^{i}_{\infty
   }=\begin{cases} {F} \mbox{ for } i=0 \\ 0 \mbox{ for } i\ne 0.\end{cases}$$
\end{theorem}
\begin{proof}
We will only prove the existence of the first spectral sequence.
The other can be done similarly. For any $\gamma$, $0 \leq \gamma
\leq n$, we write $^{i}V^{\bullet}_\gamma$ for the graded vector
spaces
\[^{i}V^{\bullet}_\gamma = \Hom^{\bullet}_{\cD}
(R_{\cE_{n}\cdots \cE_{\gamma+1}}E_{i}^{\gamma}, F)=
\Hom^{\bullet}_{\cD}(E_{i}^{ \gamma}, L_{\cE_{\gamma+1}\cdots
\cE_n}F)\] where the second equality follows from standard
properties of mutations (\cite{GK}; Pag. 12-14).

By Remark \ref{triangles}, the triangles defining the consequent
right mutations of $F$ and the consequent left mutations of $F[n]$
through $(\cE_0, \cdots, \cE_n)$ can be written as
\[ (\bigoplus_{i=1}^{\alpha_ {\gamma}} {}^{i}V_{\gamma}^{\bullet} \otimes
E_i^{\gamma})[-1] \mapright{k_{\gamma}} R_{\cE_{\gamma}\cdots
\cE_{0}} F[-1] \mapright{i_{\gamma}} R_{\cE_{\gamma-1}\cdots
\cE_{0}} F \mapright{j_{\gamma}} \bigoplus_{i=1}^{\alpha_
{\gamma}} {}^{i}V_{\gamma}^{\bullet} \otimes E_i^{\gamma}\]
\[ \bigoplus_{i=1}^{\alpha_ {\gamma}} {}^{i}V_{\gamma}^{\bullet} \otimes
E_i^{\gamma} \mapright{j^{\gamma+1}} L_{\cE_{\gamma+1}\cdots
\cE_{n}} F[n] \mapright{i^{\gamma+1}} L_{\cE_{\gamma}\cdots
\cE_{n}} F[n+1] \mapright{k^{\gamma+1}}( \bigoplus_{i=1}^{\alpha_
{\gamma}} {}^{i}V_{\gamma}^{\bullet} \otimes E_i^{\gamma})[1].
\]
We arrange them into the following big diagram:

\vskip 6mm

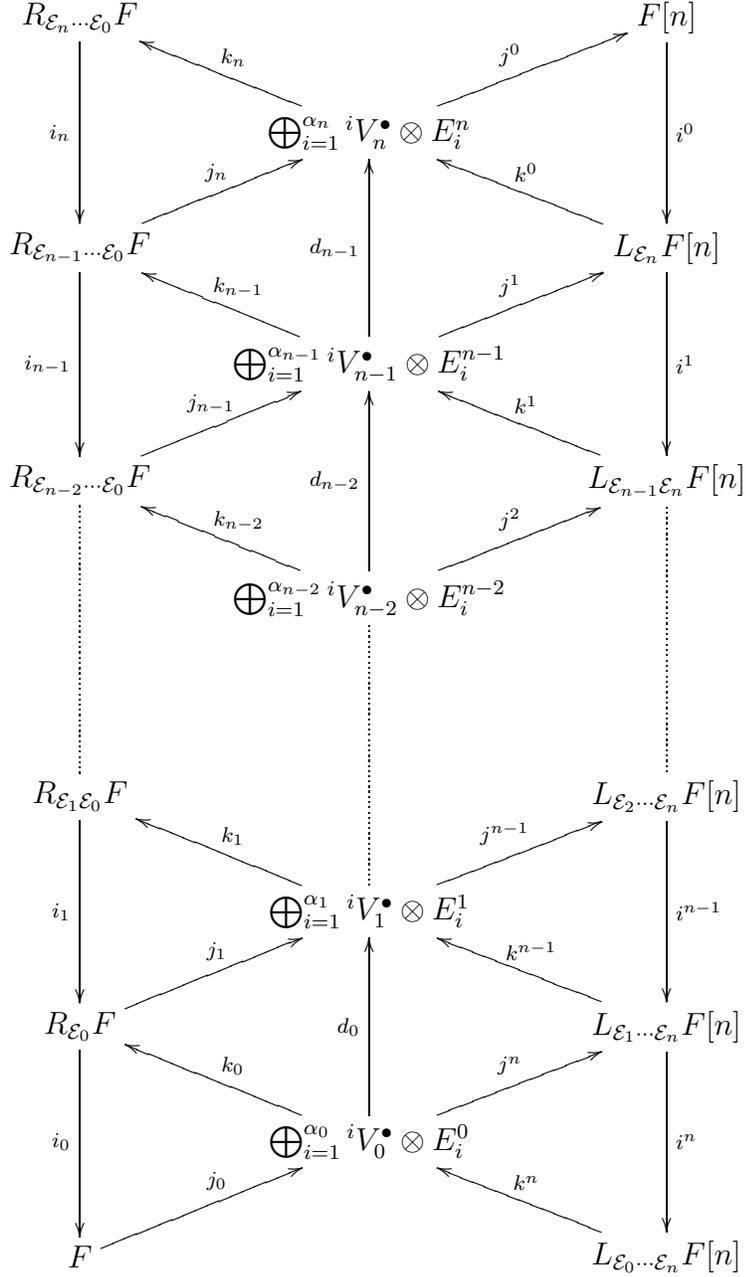
 \xymatrix{
R_{\cE_n\cdots \cE_0}F \ar[dd]_{i_n} & & F[n] \ar[dd]^{i^0}
\\
 & \ar[ul]_{k_n} \bigoplus_{i=1}^{\alpha_n}  {}^iV_n^{\bullet} \otimes E_i ^n
   \ar[ur]^{j^0} & \\
  R_{\cE_{n-1}\cdots \cE_0}F \ar[ur]^{j_n}\ar[dd]_{i_{n-1}} & &
  \ar[ul]_{k^0} \ar[dd]^{i^1} L_{\cE_n}F[n] \\
  & \ar[ul]_{k_{n-1}} \ar[uu]^{d_{n-1}}
   \bigoplus_{i=1}^{\alpha_{n-1}} {}^iV_{n-1}^{\bullet} \otimes E_i ^{n-1}  \ar[ur]^{j^1} & \\
R_{\cE_{n-2} \cdots \cE_0}F \ar[ur]^{j_{n-1}} \ar@{.}[ddd] & &
  \ar[ul]_{k^1} \ar@{.}[ddd] L_{\cE_{n-1} \cE_n}F[n] \\
  & \ar[ul]_{k_{n-2}} \ar[uu]^{d_{n-2}}
   \bigoplus_{i=1}^{\alpha_{n-2}} {}^iV_{n-2}^{\bullet} \otimes E_i ^{n-2}  \ar[ur]^{j^2} & \\
   && \\
   R_{\cE_{1} \cE_0}F  \ar[dd]_{i_1} & &
   \ar[dd]^{i^{n-1}} L_{\cE_{2} \cdots \cE_n}F[n] \\
  & \ar[ul]_{k_{1}} \ar@{.}[uuu]
   \bigoplus_{i=1}^{\alpha_{1}} {}^iV_{1}^{\bullet} \otimes E_i ^{1}  \ar[ur]^{j^{n-1}} & \\
R_{\cE_0}F \ar[ur]^{j_{1}} \ar[dd]_{i_0} & &
  \ar[ul]_{k^{n-1}} \ar[dd]^{i^n} L_{\cE_{1} \cdots \cE_n}F[n] \\
  & \ar[ul]_{k_{0}} \ar[uu]^{d_{0}}
   \bigoplus_{i=1}^{\alpha_{0}}  {}^iV_{0}^{\bullet} \otimes E_i ^{0}  \ar[ur]^{j^n} & \\
 F \ar[ur]^{j_0} & & \ar[ul]_{k^{n}} L_{\cE_{0} \cdots
 \cE_n}F[n]
 }

At this diagram, all oriented triangles along left and right
vertical borders are distinguished, the morphisms  $i_{\bullet}$
and $i^{\bullet}$ have degree one, and all triangles and rhombuses
in the central column are commutative. So, there is the following
complex, functorial on $F$,

\[L^{\bullet}: 0 \rightarrow
\bigoplus_{i=1}^{\alpha_0}{}^{i}V_0^{\bullet} \otimes E_{i}^0
 \rightarrow \bigoplus_{i=1}^{\alpha_1}{}^{i}V_1^{\bullet} \otimes E_{i}^1 \rightarrow \cdots
 \rightarrow \bigoplus_{i=1}^{\alpha_{n-1}}{}^{i}V_{n-1}^{\bullet} \otimes
 E_{i}^{n-1}
  \rightarrow \bigoplus_{i=1}^{\alpha_n}{}^{i}V_{n}^{\bullet} \otimes
  E_{i}^{n}
 \rightarrow 0   \]
 and by the above Postnikov-system we have that $F$ is a right convolution of this complex.
  Then, for an
 arbitrary linear covariant cohomological functor
 $\Phi^{\bullet}$, there exists an spectral sequence with
 $E_1$-term
 \[ _{I}E_1^{pq}=\Phi^q(L^p)\]
 situated in the square $0 \leq p,q \leq n$ and converging to
  $\Phi^{p+q}(F)$ (see \cite{Kap2}; 1.5). Since
 $\Phi^{\bullet}$ is a linear functor, we have
\begin{equation} \label{espectralgral}
 \Phi^q(L^p)= \bigoplus_{i=1}^{\alpha_p}\Phi^q(^{i}V_p^{\bullet} \otimes E_{i}^{p})=
\bigoplus_{i=1}^{\alpha_p}\bigoplus_{l}  {}^{i}V_p^l \otimes
\Phi^{q-l}(E_{i}^{p})= \bigoplus_{i=1}^{\alpha_p}
\bigoplus_{\alpha+\beta=q} {}^{i}V_p^{\alpha} \otimes
\Phi^{\beta}(E_{i}^{p}). \end{equation} In particular, if we
consider the covariant linear cohomology functor which takes a
complex to its cohomology sheaf and acts identically on pure
sheaves, i.e.
\[ \Phi^{\beta}(F)= \begin{cases} F \mbox{ for } \beta=0 \\ 0 \mbox{ for } \beta \neq 0
 \end{cases}\]
 on any pure sheaf $F$, in the square $0 \leq p,q \leq n$, we get
 \[ _{I} E_1^{pq}=\bigoplus_{i=1}^{\alpha_p}
{}^{i}V_p^{q} \otimes E_{i}^{p}= \bigoplus_{i=1}^{\alpha_p}
 \Ext^q(R_{\cE_{n} \cdots \cE_{p+1}}E _{i}^p,{F})
   \otimes E_{i}^p \]
 which converges to
 $$ _{I} E^{i}_{\infty }=
\begin{cases} {F} \mbox{ for } i=0 \\ 0 \mbox{ for } i\ne
0.\end{cases}$$ Finally, if we call $p'=p-n$, we get the spectral
sequence $$   _{I}
E_1^{p'q}=\bigoplus_{i=1}^{\alpha_{p'+n}}\Ext^q(R_{\cE_{n} \cdots
\cE_{p'+n+1}}E _{i}^{p'+n},{F})
   \otimes E_{i}^{p'+n} $$
  situated in the square $-n \le p'  \le 0$, $0 \le q \le n $
   which converges to $$ _{I} E^{i}_{\infty }=\begin{cases} {F} \mbox{ for } i=0 \\ 0 \mbox{ for } i\ne 0.\end{cases}$$
\end{proof}

\begin{remark} We want to point out that in Theorem \ref{mbe} the number of
blocks is one greater than the dimension of $X$ but a priori there
is no restriction on the length $\alpha_j$ of each block
$\cE_j=(E_1^j, \ldots,E_{\alpha_j}^j)$.
\end{remark}

\vspace{3mm} Arguing as in \cite{CMR}; Lemma 2.23 but using the
distinguished triangles (\ref{triangleleft}) and
(\ref{triangleright}) instead of exact sequences, we can prove the
following technical Lemma that will be used in next sections.

\begin{lemma}\label{tecnickey}
\label{tecnic} Let $X$ be a smooth projective variety of dimension
$n$ and let $\cB=(\cE_0, \cdots, \cE_n)$ be an $n$-block
collection of coherent sheaves on $X$. For any $i < j$ and any
locally free sheaf $F$ on $X$, it holds:

(a) $(L_{\cE_i}E_k^{j})^*=R_{\cE_i^*} {E_k^{j}}^*$ for any
$E^{j}_{k}\in \cE_j$;

(b) $(R_{\cE_j}E_k^{i})^*=L_{\cE_j^*} {E_k^{i}}^*$ for any
$E^{i}_{k}\in \cE_i$;

(c) $(R_{\cE_j}E_k^{i})\otimes F \cong R_{\cE_j \otimes F}(E_k^{i}
\otimes F)$ for any $E^{i}_{k}\in \cE_i$ and
$(L_{\cE_i}E_k^{j})\otimes F \cong L_{\cE_i \otimes F}(E_k^{j}
\otimes F)$ for any $E^{j}_{k}\in \cE_j$.
\end{lemma}


\section{Regularity with respect to $n$-blocks collections}

The goal of this section is to extend the notion of
Castelnuovo-Mumford regularity for coherent sheaves on a
projective space to coherent sheaves on an $n$-dimensional  smooth
projective variety with an $n$-block collection of coherent
sheaves on $X$ which generates $\cD$. We establish for coherent
sheaves on $\PP^n$ the agreement of the new definition of
regularity with the old one and we prove that many formal
properties of Castelnuovo-Mumford regularity continue to hold in
our more general setup.

\vskip 2mm To extend the notion of Castelnuovo-Mumford regularity,
we will first introduce the notion of helix of blocks associated
to an $m$-block collection of objects of $\cD$ as a natural
generalization of the notion of helix associated to an exceptional
collection of objects of $\cD$ introduced by J.M. Drezet and J. Le
Potier in  \cite{DLP} and L. Gorodentsev and A.N. Rudakov in
\cite{GR} (see also \cite{H}).

\begin{definition}\label{helix} Let $X$ be a smooth projective variety
and let $\cB =(\cE_0, \cE_1,\cdots ,\cE _{m})$, $\cE_j=(E_1^j,
\ldots,E_{\alpha_j}^j )$  be an $m$-block collection
  of objects of $\cD$. We extend, in both directions, the
  collection $\cB$ to an infinite sequence of blocks defining by
  induction:
  $$ \cE_{i+m}:=R^{(m)}\cE_{i-1} \mbox{ and }
  \cE_{-i}=L^{(m)}\cE_{m-i+1} \quad i>0.$$
  The collection ${\mathcal H}_{\cB }=\{ \cE_ {i} \}_{i\in \ZZ}$
is called the {\bf helix of blocks associated to $\cB$}. Each
helix ${\mathcal H}_{\cB }=\{ \cE_ {i} \}_{i\in \ZZ}$ is uniquely
recovered from any collection of $(m+1)$ of its consequent blocks
$\cB _{i}=(\cE_i, \cE_{i+1},\cdots ,\cE _{i+m})$.
\end{definition}

\begin{remark}\label{helix-helixblocks}  Let $X$ be a smooth projective
variety, let $\cB =(\cE_0, \cE_1,\cdots ,\cE _{m})$,
$\cE_j=(E_1^j, \ldots,E_{\alpha_j}^j )$,  be an $m$-block
collection
  and let ${\mathcal H}_{\cB }=\{ \cE_ {i} \}_{i\in \ZZ}$  be the
  helix of blocks associated to $\cB$. If we consider ${\mathcal H}_{\cB
  }$ just as a collection of objects of $\cD$ (we forget the
  blocks) then it turns out to be the helix associated to the
  exceptional collection \[\sigma=(\cE_0, \cE_1, \cdots, \cE_m)=(E_1^0,\cdots,E_{\alpha_0}^0, E_1^1,
 \cdots, E_{\alpha_1}^1,
\cdots, E_1^m,\cdots, E_{\alpha_m}^m). \]
\end{remark}

\begin{definition}
\label{period} Let $X$ be a smooth projective variety of dimension
$n$ with canonical bundle $K_X$. A sequence $\{\cE_i\}_{i \in
\ZZ}$ of blocks of  objects of $\cD$ will be called a {\bf helix
of blocks of period $(m+1)$} if for any $i \in \ZZ$,
\[ \cE_{i}= \cE_{i+m+1} \otimes K_X[n-m] \]
where $\cE_j\otimes K_X[n-m]$ denotes the block $(E_1^j\otimes
K_X[n-m],E_2^j\otimes K_X[n-m], \cdots ,E_{\alpha _j}^j\otimes
K_X[n-m])$  and the number in square brackets denotes the
multiplicity of the shift of an object to the left viewed as a
graded complex in $\cD$.
\end{definition}

\vspace{2mm} It follows from Remark \ref{helix-helixblocks} and
\cite{Bo}; Theorem 4.1 that if $X$ is a smooth projective variety
of dimension $n$ and $\cB =(\cE_0, \cE_1,\cdots ,\cE _{m})$,
$\cE_j=(E_1^j, \ldots,E_{\alpha_j}^j )$,  is an $m$-block
collection of objects  of $\cD$, then the helix of blocks
${\mathcal H}_{\cB }=\{ \cE_ {i} \}_{i\in \ZZ}$   associated to
$\cB$ is an helix of blocks of period $m+1$, i.e., for any $i\in
\ZZ$, $\cE_i=\cE_{i+m+1}\otimes K_X[n-m]$.

 As an immediate consequence of Proposition
\ref{bonespropietats}, we have:

\vspace{2mm}

\begin{corollary}
\label{excellent2} Let $X$ be a smooth projective variety of
dimension $n$, let $\cB=(\cE_0, \cdots, \cE_n)$, $\cE_j=(E_1^j,
\ldots,E_{\alpha_j}^j )$  be an $n$-block collection of coherent
sheaves on $X$ which generates $\cD$
 and denote by
$\cH_{\cB}=\{ \cE_ {i} \}_{i\in \ZZ}$ the helix of blocks
associated to $\cB$. Then, $\cH_{\cB}$ is an helix of period $n+1$
and any $n$-block collection $\cB _i=(\cE _{i},\cE _{i+1},\cdots
,\cE _{i+n})$ of $n+1$ subsequent blocks is an  $n$-block
collection of coherent sheaves on $X$ which generates $\cD$.
\end{corollary}

\vskip 2mm Let $X$ be a smooth projective variety of dimension $n$
and let $\cB =(\cE_0, \cE_1,\cdots ,\cE _{n})$, $\cE_j=(E_1^j,
\ldots,E_{\alpha_j}^j )$  be an $n$-block collection
  of coherent sheaves on $X$ which generates $\cD$.
  Associated to $\cB $ we have a  helix
of blocks ${\mathcal H}_{\cB }=\{ \cE_ {i} \}_{i\in \ZZ}$; and for
any $n$-block collection $\cB _i=(\cE _{i},\cE _{i+1},\cdots ,\cE
_{i+n})$ of $n+1$ subsequent blocks and any coherent
${\cO}_{X}$-module $F$ we have a spectral sequence (See Theorem
\ref{mbe})

\begin{equation}\label{sucesiopectral}
   _{I}^{i} E_1^{pq}=
\begin{cases}
   \bigoplus_{s=1}^{\alpha_{p+n+i}}\Ext^q(R_{\cE_{i+n} \cdots \cE_{i+p+n+1}}E _{s}^{i+p+n},F)
   \otimes E_{s}^{i+p+n}   & \mbox{if} \quad -n \leq p \leq -1 \\
\bigoplus_{s=1}^{\alpha_{i+n}}\Ext^q(E _{s}^{i+n},F)
   \otimes E_{s}^{i+n}   & \mbox{if} \quad p=0 \end{cases}
   \end{equation}
 situated in the square $0\le q\le n $, $-n
\le p \le 0$ which converges to $$ E^{r}_{\infty }=
\begin{cases} {F} \mbox{ for } r=0 \\ 0 \mbox{ for } r\ne 0.\end{cases}$$

\begin{definition} \label{new} Let $X$ be a smooth projective  variety of dimension $n$
 with an  $n$-block collection $\cB =(\cE_0, \cE_1,\cdots ,\cE _{n})$,
  $\cE_i=(E_1^i, \ldots,E_{\alpha_i}^i )$ of coherent sheaves on $X$ which generates
  $\cD$, let
  ${\mathcal H}_{\cB }=\{ \cE_ {i} \}_{i\in \ZZ}$ be the helix of
  blocks associated to $\cB$    and
 let $F$ be a coherent ${\cO}_{X}$-module. We say that
 $F$ is {\bf $m$-regular with respect to $\cB$} if for $q>0$
 we have
$$ \begin{cases}
   \bigoplus_{s=1}^{\alpha_{-m+p}}\Ext^q(R_{\cE_{-m} \cdots \cE_{-m+p+1}}E _{s}^{-m+p},F)
   =0  & \mbox{for} \quad -n \leq p \leq -1 \\
\bigoplus_{s=1}^{\alpha_{-m}}\Ext^q(E _{s}^{-m},F)
   =0   & \mbox{for} \quad p=0. \end{cases}$$
\end{definition}

\vskip 2mm So, $F$ is $m$-regular with respect to $\cB $ if
$^{-n-m}E_1^{pq}=0$ for $q>0$ in (\ref{sucesiopectral}). In
particular, if $F$ is $m$-regular with respect to $\cB$ the
spectral sequence $^{-n-m}E_1^{pq}$ collapses at $E_2$  and we get
the following exact sequence:

\begin{equation} \label{re}
0 \longrightarrow {\cL}_{-n} \longrightarrow \cdots
\longrightarrow {\cL}_{-1} \longrightarrow {\cL}_0 \longrightarrow
F \longrightarrow 0
\end{equation}
where $${\cL}_p=
\begin{cases}
   \bigoplus_{s=1}^{\alpha_{-m+p}}H^0(X,(R_{\cE_{-m} \cdots \cE_{-m+p+1}}E _{s}^{-m+p})^*\otimes
   F)
   \otimes E_{s}^{-m+p}   & \mbox{if} \quad -n \leq p \leq -1 \\
\bigoplus_{s=1}^{\alpha_{-m}}H^0((E _{s}^{-m})^*\otimes F)
   \otimes E_{s}^{-m}   & \mbox{if} \quad p=0. \end{cases}$$

\vspace{3mm}

\begin{definition}\label{B-regular}  Let $X$ be a smooth
projective variety of dimension $n$
 with an $n$-block collection $\cB =(\cE_0, \cE_1,\cdots ,\cE _{n})$,
  $\cE_j=(E_1^j, \ldots,E_{\alpha_j}^j )$ of coherent sheaves on $X$ which generates $\cD$
   and
 let $F$ be a coherent ${\cO}_{X}$-module. We define the
 {\bf regularity of $F$ with respect to $\cB$} (or {\bf $\cB$-regularity of $F$}), $Reg_{\cB}(F)$, as the least integer $m$ such
 that $F$ is $m$-regular with respect to $\cB$ in the sense of Definition
 \ref{new}. We set $Reg_{\cB}(F)=-\infty $ if there is no such
integer.
\end{definition}

\begin{remark} Let $X$ be a smooth
projective variety of dimension $n$. Since a geometric collection
 $\sigma =(E_0, E_1,\cdots ,E _{n})$ is an $n$-block collection of type (1, $\cdots $,1),
 Definitions
 \ref{new} and \ref{B-regular} extend the definition of regularity with respect to a geometric collection
 introduced by the authors in \cite{CMR}.
\end{remark}

\begin{example} We consider the $n$-block collection $\cB=(\cO_{\PP^n},
\cO_{\PP^n}(1), \cdots ,\cO_{\PP^n}(n))$ on $\PP^n$ and the
associated  helix ${\mathcal H}_{\cB }=\{\cO_{\PP^n}(i)\}_{i\in
\ZZ}$. According to Example \ref{exampledualbase} (1), the left
dual $n$-block collection of an $n$-block collection
$\cB_i=(\cO_{\PP^n}(i), \cO_{\PP^n}(i+1), \cdots
,\cO_{\PP^n}(i+n))$ of $n+1$ subsequent blocks of ${\mathcal
H}_{\cB }$ is
$$(\cO_{\PP^n}(i+n), R^{(1)}\cO_{\PP^n}(i+n-1),\cdots
,R^{(j)}\cO_{\PP^n}(i+n-j),\cdots ,R^{(n)}\cO_{\PP^n}(i))$$
$$=(\cO_{\PP^n}(i+n), T_{\PP^n}(i+n-1),\cdots
,\wedge^{j}T_{\PP^n}(i+n-j),\cdots ,\wedge^{n}T_{\PP^n}(i)).$$
Therefore, for any coherent sheaf $F$ on $\PP^n$ our definition
reduces to say: $F$ is $m$-regular with respect to $\cB$ if
$\Ext^q(\wedge^{-p}T(-m+p),F)=H^{q}(\PP^n,\Omega ^{-p}(m-p)\otimes
F)=0$ for all $q>0$ and all $p$, $-n\le p\le 0$.
\end{example}

We will now compute the regularity with respect to $\cB =(\cE_0,
\cE_1,\cdots ,\cE _{n})$ of any coherent sheaf $E _{t}^{i}\in
\cE_{i}$.

\begin{proposition} \label{regEi} Let $X$ be a smooth projective  variety of dimension $n$
 with an $n$-block collection $\cB=(\cE_0, \cE_1,\cdots ,\cE _{n})$, $\cE_j=(E_1^j,
\ldots,E_{\alpha_j}^j )$,   which generates $\cD$ and
 let ${\mathcal H}_{\cB}=\{ \cE _{i}\} _{i\in \ZZ}$ be the
 associated
 helix. Then, for any $i\in \ZZ$ and  any $E_{t}^{i}\in \cE_{i}$, $Reg_{\cB }(E _{t}^{i})=-i$.
\end{proposition}
\begin{proof} First of all we will see that $Reg_{\cB}(E_t^{i})\le
-i$. By the orthogonality relation (\ref{orto1r}), for $q>0$ and
$-n\le p \le 0$, we have
$$\bigoplus _{s=1}^{\alpha_{i+p}}\Ext ^q(R_{\cE_{i}\cdots
\cE_{i+p+1}}E_s^{i+p},E_t^{i})=0 \mbox{ for } -n\le p \le -1,
\mbox{ and }$$
$$\bigoplus _{s=1}^{\alpha_{i}}\Ext ^q(E_s^{i},E_t^{i})=0 \mbox{ for }  p=0.$$
So, $E_t^{i}$ is $(-i)$-regular with respect to $\cB$ or,
equivalently, $Reg_{\cB} (E_t^{i})\le -i$.

Let us now see that $E_t^{i}$ is not $(-i-1)$-regular with respect
to $\cB$. To this end, it is enough to see that
$$\bigoplus_{s=1}^{\alpha_i}\Ext^1(R_{\cE_{i+1}}E_s^{i},E_t^{i})\ne 0.$$

To prove it, we write  $i=\alpha n+j$ with $0\le j <n$, $\alpha
\in \ZZ$. We consider the $n$-block collection $$\cB_{\alpha
n}=(\cE _{\alpha n},\cE _{\alpha n+1},\cdots ,\cE _{\alpha
n+j}=\cE _{i},\cdots ,\cE _{\alpha n+n})$$ and  its left dual
$n$-block collection
$$(\cE _{\alpha n+n},R^{(1)}\cE _{\alpha n+n-1},\cdots
,R^{(n-j)}\cE _{\alpha n+j},\cdots ,R^{(n)}\cE _{\alpha n}).$$ By
Remark \ref{existence-dual} and the equality (\ref{orto2})

 $$\Ext^{n-j}(R^{(n-j)}E_t^{i}
, E_t^{i} )= \Ext^{n-j}(R_{\cE _{\alpha n+n}\cdots \cE
_{i+2}\cE_{i+1}}E _t^{i}, E _{t}^{i})=\CC.$$

So, if $j=n-1$ we are done. Assume $j<n-1$. Since $\Ext^q(E
_{s}^{i},E_{t}^{i})= 0$ for $q>0$ and $1\le s\le \alpha_i$,
applying the contravariant functor $\Hom (.,E _{t}^{i})$ to the
exact sequence
$$0 \longrightarrow  R_{\cE _{\alpha n+n-1} \cdots \cE
_{i+1}} E _{t}^{i} \longrightarrow
\bigoplus_{s=1}^{\alpha_i}\Hom^*(R_{\cE _{\alpha n+n-1} \cdots \cE
_{i+1}} E _{t}^{i}, E _{s}^{i})\otimes E _{s}^{i} \longrightarrow
R^{(n-j)}E _{t}^{i} \longrightarrow 0$$ we obtain
$$\Ext^{n-j-1}(R_{\cE _{\alpha n+n-1} \cdots \cE _{i+1}} E _{t}^{i},E_{t}^{i})=\CC.$$

We repeat the process using the consequent right mutations and we
get
$$\Ext^{n-j-k}(R_{\cE _{\alpha n+n-k} \cdots \cE _{i+1}} E
_{t}^{i},E_{t}^{i})=\CC$$ for $0\le k\le n-1-j$. In particular,
$$\Ext^1(R_{\cE _{i+1}}E _t^{i},E_t^{i})=\CC $$ and, hence, $$
\bigoplus_{s=1}^{\alpha_i}\Ext^1(R_{\cE_{i+1}}E_s^{i},E_t^{i})\ne
0
$$ which implies that
$E _t^{i}$ is not $(-i-1)$-regular and we conclude that $Reg_{\cB
}(E _t^{i})=-i$.
\end{proof}

Let us now compare our new definition of regularity with the
previous ones.

\vskip 3mm \noindent {\bf Castelnuovo-Mumford regularity}. In
\cite{M}, Lecture 14, D. Mumford defined the notion of regularity
for a coherent sheaf over a projective space. Let us recall it.

\vspace{3mm}

\begin{definition} \label{defCM} A coherent sheaf $F $ on $\PP^n$ is said to
be {\bf $m$-regular in the sense of Castelnuovo-Mumford} if
$H^{i}(\PP^n,F (m-i))=0$ for $i>0$. We define the
Castelnuovo-Mumford regularity of $F$, $Reg^{CM}(F)$, as  the
least integer $m$ such that $F$ is $m$-regular. We say that the
Castelnuovo-Mumford regularity is $-\infty $ if  such integer does
not exist.
\end{definition}

\vspace{3mm}

Let us now establish for coherent sheaves on $\PP^n$ the agreement
of the $\cB$-regularity in the sense of Definition \ref{new} with
Castelnuovo-Mumford definition.

\begin{proposition} \label{CM}  A coherent sheaf $F $ on $\PP^n$ is
 $m$-regular in the sense of Castelnuovo-Mumford if and only if it is
 $m$-regular with respect to the $n$-block collection
$\cB =(\cO_{\PP^n}, $ $\cO_{\PP^n}(1),\cdots ,\cO_{\PP^n}(n))$
 in the sense of Definition \ref{new}. Hence, we have
 $$ Reg_{\cB }(F)=Reg^{CM}(F).$$
\end{proposition}

\begin{proof} Since $\cB$ is not only an $n$-block collection but
also a geometric collection, the result follows from \cite{CMR};
Proposition 4.6.
\end{proof}

\vskip 3mm \noindent {\bf Chipalkatti's regularity}. In \cite{Ch},
Definition 1.1, J.V. Chipalkatti introduced  the notion of
regularity for a coherent sheaf on a Grassmannian variety and he
showed that when the Grassmannian is a projective space his
definition of regularity agrees with Castelnuovo-Mumford
regularity. We will now see that the notion of regularity
introduced in Definition \ref{B-regular} is closely related to
Chipalkatti's regularity but they do  not coincide. Let us recall
Chipalkatti' definition.

\begin{definition} \label{chipa} A sheaf $F$ on $X=Gr(k,n)$ is regular if
$H^q(X,F\otimes \Sigma ^{\beta}Q^*)=0$ for all $\beta $ such that
$k\ge \beta_1\ge \cdots \beta_{n-k}\ge 0$ and all $q\ge 1$. It is
said to be $m$-regular if $F\otimes \cO_X(m)$ is regular.
$Reg^C(F)$ is the least integer $m$ such that $F(m)$ is regular,
set $Reg^C(F)=-\infty $ if there is no such integer.
\end{definition}

Proposition \ref{CM} together with \cite{Ch}; Theorem 1.5 and 1.6,
establishes, when the Grassmannian is a projective space, the
agreement of Chipalkatti's definition (Definition \ref{chipa}),
Castelnuovo-Mumford's definition (Definition \ref{defCM}) and our
definition (Definition \ref{new}). Nevertheless, next example
shows that, in general, Chipalkatti's definition and our
definition do not coincide.

\begin{example} We consider the 4-dimensional Grassmann variety
$X=Gr(2,4)$ and the 4-block collection described in Example
\ref{exempleblock} (2) tensored with $\cO_X(2)$
$$\cB=(\cE_0,\cE_1,\cE_2,\cE_3,\cE_4)=(\Sigma^{(2,2)}\cS\otimes \cO_X(2), \Sigma^{(2,1)}\cS\otimes
\cO_X(2),$$$$ (\Sigma^{(2,0)}\cS\otimes \cO_X(2),
\Sigma^{(1,1)}\cS\otimes \cO_X(2)), \Sigma^{(1,0)}\cS\otimes
\cO_X(2), \Sigma^{(0,0)}\cS\otimes \cO_X(2)).$$ Since
$\Sigma^{(2,2)}\cS\otimes \cO_X(2)=\cO_X$ and
$\Sigma^{(0,0)}\cS\otimes \cO_X(2)=\cO_X(2)$, by \cite{Ch};
Example 1.3(a), $Reg^C(\Sigma^{(0,0)}\cS\otimes \cO_X(2))=-2$ and
$Reg^C(\Sigma^{(2,2)}\cS\otimes \cO_X(2))=0$. On the other hand,
by Proposition \ref{regEi}, $Reg_{\cB}(\Sigma^{(0,0)}\cS\otimes
\cO_X(2))=-4$ and $Reg_{\cB}(\Sigma^{(2,2)}\cS\otimes
\cO_X(2))=0$.

\end{example}

\vskip 3mm \noindent {\bf Hoffman-Wang regularity}. In \cite{HW},
J.W. Hoffman and H.H. Wang  introduced a multigraded variant of
the Castelnuovo-Mumford regularity and we will devote section 5 of
this paper to relate it to our new definition of regularity.

\vskip 4mm To emphasize the similarities between the new notion of
regularity and the original definition in chapter 14 of \cite{M},
we will end this section proving that the basic formal properties
of Castelnuovo-Mumford regularity of coherent sheaves over
projective spaces remains to be true in this new setting.

\vspace{3mm}

\begin{proposition}\label{sup}  Let $X$ be a smooth projective  variety of
dimension $n$ with an $n$-block collection  of coherent sheaves
$\cB=(\cE_0, \cE_1,\cdots ,\cE _{n})$, $\cE_j=(E_1^j,
\ldots,E_{\alpha_j}^j )$,
 which generates $\cD$
  and
 let $F$ be a coherent ${\cO}_{X}$-module. If $F $ is
$m$-regular with respect to $\cB $ then  the canonical map
$\bigoplus_{s=1}^{\alpha_{-m}}\Hom(E _{s}^{-m}, F)
   \otimes E_{s}^{-m} \twoheadrightarrow F$ is
surjective and $F$ is $k$-regular with respect to $\cB $ for any
$k\ge m$ as well.
\end{proposition}

\begin{proof}
The first assertion follows from the exact sequence (\ref{re}). To
prove the second assertion it is enough to check it for $k=m+1$.
Since $F$ is  $m$-regular with respect to $\cB$ we have for $q>0$
\begin{equation}\label{m-regularby hypothesis} \begin{cases}
   \bigoplus_{s=1}^{\alpha_{-m+p}}\Ext^q(R_{\cE_{-m} \cdots \cE_{-m+p+1}}E _{s}^{-m+p},F)
   =0  & \mbox{if} \quad -n \leq p \leq -1 \\
\bigoplus_{s=1}^{\alpha_{-m}}\Ext^q(E _{s}^{-m},F)
   =0   & \mbox{if} \quad p=0. \end{cases} \end{equation}
   In order to see that $F$ is $(m+1)$-regular with respect to $\cB$ we
   have to prove
   \begin{equation}\label{m+1-regular} \begin{cases}
   \bigoplus_{s=1}^{\alpha_{-m-1+p}}\Ext^q(R_{\cE_{-m-1} \cdots \cE_{-m+p}}E _{s}^{-m-1+p},F)
   =0  & \mbox{for }  q>0, \mbox{ } -n \leq p \leq -1 \\
\bigoplus_{s=1}^{\alpha_{-m-1}}\Ext^q(E _{s}^{-m-1},F)
   =0   & \mbox{for }  q>0, \mbox{  } p=0. \end{cases}
   \end{equation}
   Using the equalities (\ref{m-regularby hypothesis}) and applying, for any $s$, $1\le s \le \alpha_{-m-1}$,
the functor $\Hom(.,F)$ to the exact sequence
$$ 0\longrightarrow E _{s}^{-m-1} \longrightarrow \oplus_{t=1}^{\alpha_{-m}}\Hom^*(E _{s}^{-m-1}
,E _{t}^{-m})\otimes E _{t}^{-m}\longrightarrow R_{\cE
_{-m}}E_s^{-m-1}\longrightarrow  0$$ we obtain
$$\Ext^q(E _{s}^{-m-1},F)
   =0   \mbox{ for }  q>0 \mbox{ and } 1\le s \le \alpha_{-m-1} $$
   and,
   hence, $$\bigoplus_{s=1}^{\alpha_{-m-1}}\Ext^q(E _{s}^{-m-1},F)
   =0   \mbox{ for }  q>0 .$$

   Using again the equalities (\ref{m-regularby hypothesis}) and,
   for any $s$, $1\le s \le \alpha_{-m-2}$,
   the exact sequence
   $$ 0\longrightarrow R_{\cE _{-m-1}}E _{s}^{-m-2} \longrightarrow
\oplus_{t=1}^{\alpha_{-m}}\Hom^*(R_{\cE _{-m-1}}E _{s}^{-m-2},E_t
^{-m})\otimes E _t^{-m}\longrightarrow R_{\cE _{-m}\cE _{-m-1}}E
_{s}^{-m-2}\longrightarrow 0$$ we get
$$\bigoplus_{s=1}^{\alpha_{-m-2}}\Ext^q(R_{\cE_{-m-1}} E _{s}^{-m-2},F)
   =0   \mbox{ for any }  q>0.$$

Going on and using the consequent right  mutations  of blocks, we
get for all $p$, $-n+1 \leq p \leq -1,$
$$\bigoplus_{s=1}^{\alpha_{-m-1+p}}\Ext^q(R_{\cE_{-m-1} \cdots
\cE_{-m+p}}E _{s}^{-m-1+p},F)
   =0 \mbox{ for all } q>0.  $$
Therefore, it only remains to see that
$$\bigoplus_{s=1}^{\alpha_{-m-1-n}}\Ext^q(R_{\cE_{-m-1} \cdots
\cE_{-m-n}}E _{s}^{-m-1-n},F)
   =0 \mbox{ for all } q>0. $$

The vanishing of these last $\Ext$'s groups follows again from the
equalities (\ref{m-regularby hypothesis})  taking into account
that, by Definition \ref{helix} \[ R_{\cE _{-m-1}}R_{\cE
_{-m-2}}\cdots
 R_{\cE_{-m-n}}\cE_{-m-n-1}= R^{(n)} \cE_{-m-n-1}=\cE _{-m}. \]
\end{proof}

\vspace{3mm}

\begin{proposition}
\label{propietats} Let $X$ be a smooth projective  variety of
dimension $n$ with an $n$-block collection  of coherent sheaves
$\cB=(\cE_0, \cE_1,\cdots ,\cE _{n})$, $\cE_j=(E_1^j,
\ldots,E_{\alpha_j}^j )$,
 which generates $\cD$. Let $F$
 and $G$ be coherent ${\cO}_{X}$-modules and
 let
\begin{equation} \label{regsuc} 0  \longrightarrow F_1 \longrightarrow F_2 \longrightarrow F_3
\longrightarrow 0\end{equation} be an exact sequence of coherent
${\cO}_{X}$-modules. Then,

 \begin{itemize}
 \item[(a)] $Reg_{\cB }(F_2) \leq max\{Reg_{\cB }(F_1),Reg_{\cB
 }(F_3)\}$,
  \item[(b)] $Reg_{\cB }(F \oplus G)=max\{Reg_{\cB }(F),Reg_{\cB
  }(G)\}$.
 \end{itemize}
\end{proposition}
\begin{proof}
(a)  Let $m= max\{Reg_{\cB }(F_1),Reg_{\cB
 }(F_3)\}$. Since, by Proposition \ref{sup}, $F_1 $ and $F_3 $ are
both $m$-regular with respect to $\cB $ considering the long exact
sequences
\[ \cdots \longrightarrow  \oplus_{s=1}^{\alpha_{-m+p}}\Ext^q(R_{\cE_{-m} \cdots \cE_{-m+p+1}}E _{s}^{-m+p},F_1) \longrightarrow
  \oplus_{s=1}^{\alpha_{-m+p}}\Ext^q(R_{\cE_{-m} \cdots \cE_{-m+p+1}}E
  _{s}^{-m+p},F_2)\]
 \[  \longrightarrow
 \bigoplus_{s=1}^{\alpha_{-m+p}}\Ext^q(R_{\cE_{-m} \cdots \cE_{-m+p+1}}E _{s}^{-m+p},F_3)
 \longrightarrow \cdots \] and
 \[ \cdots \longrightarrow  \oplus_{s=1}^{\alpha_{-m}}\Ext^q(E _{s}^{-m},F_1) \longrightarrow
  \oplus_{s=1}^{\alpha_{-m}}\Ext^q(E _{s}^{-m},F_2)\longrightarrow
 \oplus_{s=1}^{\alpha_{-m}}\Ext^q(E _{s}^{-m},F_3)
 \longrightarrow \cdots \]
associated to (\ref{regsuc}) we get

$$ \begin{cases}
   \bigoplus_{s=1}^{\alpha_{-m+p}}\Ext^q(R_{\cE_{-m} \cdots \cE_{-m+p+1}}E _{s}^{-m+p},F_2)
   =0  & \mbox{if} \quad -n \leq p \leq -1 \\
\bigoplus_{s=1}^{\alpha_{-m}}\Ext^q(E _{s}^{-m},F_2)
   =0   & \mbox{if} \quad p=0, \end{cases} $$
which implies that $Reg_{\cB}(F_2)\le m$.

 (b) It easily follows from the additivity of the functor
$\Ext^q(R_{\cE_{-m} \cdots \cE_{-m+p+1}}E _{s}^{-m+p},.)$.
\end{proof}

\section{Regularity of sheaves on multiprojective spaces}

In this section, we will restrict our attention to coherent
sheaves over multiprojective spaces $X=\PP^{n_1}\times \cdots
\times \PP^{n_r}$ and we will relate our definition of regularity
to the multigraded variant of the Castelnuovo-Mumford regularity
introduced by Hoffman and Wang \cite{HW} (See also \cite{MS}).

\vskip 2mm We first fix the notation we need in this section. For
each integer $i>0$, let $$St_i=\{(l,s)\in \ZZ^2 \mid l+s=-1-i,
l<0, s<0\}$$
$$=\{(-i,-1),(-i+1,-2),\cdots ,(-2,-i+1),(-1,-i)\},$$
for $i\le 0$, let $$St_i=\{(l,s)\in \ZZ^2 \mid l+s=-i, l\ge 0,
s\ge 0 \}$$
$$=\{(-i,0),(-i-1,1),\cdots ,(1,-i-1),(0,-i)\}.$$ For each $(p,p')\in \ZZ^2$, let
$St_i(p,p')=(p,p')+St_i$.

\begin{definition}\label{hof-wang} Let $F$ be  a coherent sheaf on $X=\PP^m\times
\PP^n$. We say that $F$ is $(p,p')$-regular if, for all $i\ge 1$,
$$H^{i}(X,F(k,k'))=0$$ whenever $(k,k')\in St_{i}(p,p')$.
\end{definition}

\begin{remark} Definition \ref{hof-wang} generalizes in an obvious
way to coherent sheaves on multiprojective spaces $\PP^{n_1}\times
\cdots \times \PP^{n_r}$.
\end{remark}

Set $d=m+n$, $X=\PP^m \times \PP^n$ and denote by $\cB =(\cE_0,
\cE_1,\cdots ,\cE _{d})$ the $d$-block collection where for any
$0\le j \le d$, we have $$\cE_j=\{ \cO_X(a,b) \mid a+b=j-d, 0\ge a
\ge -m, 0 \ge b \ge -n \}$$ and we set $\alpha_j:= \sharp \cE_j$.
By Remark \ref{existence-dual}, the left dual $d$-block collection
$$(R^{(0)}\cE_d,R^{(1)}\cE_{d-1}, \cdots , R^{(j)}\cE_{d-j},
\cdots ,R^{(d)}\cE_0)$$ of $\cB $ is univocally determined by the
orthogonality relations (\ref{orto1}) and (\ref{orto2}), and an
intricate computation using K\"{u}nneth formula for locally free
sheaves on algebraic varieties shows that for any $\cO_X(a,b)\in
\cE_{d-j}$ and any $0\le j\le d$ we have
\begin{equation}\label{dualpnpm} R^{(j)}\cO_X(a,b)=\wedge^{-a}T_{\PP^m}(a)\boxtimes
\wedge^{-b}T_{\PP^n}(b).\end{equation}

\begin{lemma} \label{lematecnico} With the above notation, let
  ${\mathcal H}_{\cB }=\{ \cE_ {i} \}_{i\in \ZZ}$ be the helix of
  blocks associated to $\cB$. Let us denote by $\cB_{k(d+1)}$ the $d$-block collection
  of $d+1$ consecutive blocks $(\cE_{k(d+1)}, \cE_{k(d+1)+1},\cdots ,\cE
  _{k(d+1)+
  d})$. Then,
  we have

  (1) $\cE_{k(d+1)+i}=
   \{ \cO_X(a+k(m+1),b+k(n+1))
  \mid a+b=i-d, 0\ge a \ge -m, 0 \ge b \ge -n
\}$

 (2) The left dual $d$-block collection of  $\cB_{k(d+1)}$ is
$$(R^{(0)}\cE_{k(d+1)+d},R^{(1)}\cE_{k(d+1)+d-1}, \cdots ,
R^{(j)}\cE_{k(d+1)+d-j}, \cdots R^{(d)}\cE_{k(d+1)})$$ where for
any $\cO_X(a+k(m+1),b+k(n+1))\in \cE_{k(d+1)+d-j}$
$$R^{(j)}\cO_X(a+k(m+1),b+k(n+1))=\wedge^{-a}T_{\PP^m}(a+k(m+1))\boxtimes
\wedge^{-b}T_{\PP^n}(b+k(n+1)).$$

\end{lemma}

\begin{proof} (1) Applying Corollary \ref{excellent2}, we get $\cE_{k (d+1)+i}=
\cE_{i}\otimes K_X ^{-k }$ and the result follows taking into
account that $K_X
 ^{-\lambda}=\cO_X(\lambda(m+1),\lambda(n+1))$.

 (2)  Straightforward computation taking into account that the left
 dual $d$-block collection of $\cB_0=(\cE_{0}, \cE_{1},\cdots ,\cE
  _{d})$ is determined (up to isomorphism) by (\ref{dualpnpm}),
  Lemma \ref{tecnickey}
  and the equalities
  \[ \begin{array}{rl} R^{(j)}\cE_{k(d+1)+d-j}:= &
R_{\cE_{k(d+1)+d}\cdots \cE_{k(d+1)+d-j+1}}\cE_{k(d+1)+d-j} \\ = &
 R_{\cE_{d}\otimes
K_X^{-k}\cdots \cE_{d-j+1}\otimes K_X^{-k}}\cE_{d-j}\otimes
K_X^{-k} \\ = & R^{(j)}(\cE_{d-j}\otimes K_X^{-k})\\
= &
 (R_{\cE_{d}
\cdots \cE_{d-j+1}}\cE_{d-j})\otimes K_X^{-k} \\ = &
  (R^{(j)}\cE_{d-j})\otimes K_X^{-k}.\end{array} \]
\end{proof}



We have the following technical lemma.

\begin{lemma}\label{key1}  Let $F$ be a coherent sheaf on $X=\PP^m\times
\PP^n$. We have:
\begin{itemize} \item[(a)]  $F$ is $(p,p')$-regular in
the sense of Hoffman and Wang if and only if $F(p,p')$ is
$(0,0)$-regular in the sense of Hoffman and Wang. \item[(b)] $F$
is $k(d+1)+t$-regular with respect to $\cB$ in the sense of
Definition \ref{new} if and only if $F(k(m+1),k(n+1))$ is
$t$-regular with respect to $\cB$ in the sense of Definition
\ref{new}.
\end{itemize}
\end{lemma}

\begin{proof} (a) It obviously follows from Definition
\ref{hof-wang}.

(b) By Definition \ref{new}, $F$ is  $k(d+1)+t$-regular with
respect to $\cB$ if and only if for $q>0$ we have \[ \begin{cases}
   \bigoplus_{s=1}^{\beta_{-k(d+1)-t+p}}\Ext^q(R_{\cE_{-k(d+1)-t} \cdots
   \cE_{-k(d+1)-t+p+1}}E _{s}^{-k(d+1)-t+p},F)
   =0  & \mbox{if} \quad -d \leq p \leq -1 \\
\bigoplus_{s=1}^{\beta_{-k(d+1)-t}}\Ext^q(E _{s}^{-k(d+1)-t},F)
   =0   & \mbox{if} \quad p=0. \end{cases}\]

Since $\cE_{-k(d+1)+i}\cong \cE_{i}\otimes K_X^k$, applying Lemma
\ref{tecnickey}, it is equivalent to say that for $q>0$, we have

\[ \begin{cases}
   \bigoplus_{s=1}^{\beta_{-t+p}}\Ext^q(R_{\cE_{-t} \cdots \cE_{-t+p+1}}E _{s}^{-t+p}\otimes K_X^{k},F)
   =0  & \mbox{if} \quad -d \leq p \leq -1 \\
\bigoplus_{s=1}^{\beta_{-t}}\Ext^q(E_{s}^{-t}\otimes K_X^{k},F)
   =0   & \mbox{if} \quad p=0 \end{cases}\]
or,  equivalent,  for $q>0$, we have
 \[ \begin{cases}
   \bigoplus_{s=1}^{\beta_{-t+p}}\Ext^q(R_{\cE_{-t} \cdots \cE_{-t+p+1}}E _{s}^{-t+p},F\otimes K_X^{-k})
   =0  & \mbox{if} \quad -d \leq p \leq -1 \\
\bigoplus_{s=1}^{\beta_{-t}}\Ext^q(E_{s}^{-t},F\otimes K_X^{-k})
   =0   & \mbox{if} \quad p=0. \end{cases}\]
   which means that $F(k(m+1),k(n+1))$ is $t$-regular with
respect to $\cB$ in the sense of Definition \ref{new}.
\end{proof}

 We are now ready to state the main result of
this section.

\begin{theorem}\label{mainthm2} Let $F$ be a coherent sheaf on $X=\PP^m\times
\PP^n$ and set $d=n+m$. Then $F$ is $(0,0)$-regular in the sense
of Hoffman and Wang if and only if $F$ is $(-d)$-regular with
respect to $\cB$ in the sense of Definition \ref{new}.
\end{theorem}

\begin{proof} According to Definitions  \ref{hof-wang} and
\ref{new} we have to see that
\begin{equation}\label{HW1} H^{i}(X,F(r,s))=0 \mbox{ for all }
i>0, r+s=-i-1, r< 0, s< 0
\end{equation}
if and only if for $q>0$ we have
\begin{equation}\label{1CMR}
 \begin{cases}
   \oplus_{s=1}^{\alpha_{p+d}}\Ext^q(R_{\cE_{d} \cdots
   \cE_{p+1+d}}E_{s}^{p+d},F)
   =\oplus_{a+b=p \atop { 0\ge a\ge -m \atop 0\ge b\ge -n} }H^q( \Omega^{-a}_{\PP^m}(-a)
   \boxtimes \Omega^{-b}_{\PP^n}(-b)\otimes F)
   =0 \\ \hspace{10.3cm} \mbox{ for }  -d \leq p \leq -1; \\
\oplus_{s=1}^{\alpha_{d}}\Ext^q(E_{s}^{d},F)=H^q(X,F)
   =0    \mbox{ for }  p=0 . \end{cases}
\end{equation}

Let us first see that (\ref{HW1}) implies (\ref{1CMR}). Since by
\cite{HW}; Proposition 2.7, any $(0,0)$-regular sheaf is also
$(p,p')$-regular for $p\ge 0$, $p'\ge 0$, we have
\begin{equation} \label{HW2}
H^{i}(X,F(r,s))=0 \mbox{ for all } i>0, r+s\ge -i-1, r\ge -i, s\ge
-i. \end{equation}

We will see that (\ref{HW2}) implies the following stronger result
\begin{equation}
\label{2CMR} H^q( \Omega^{-a}_{\PP^m}(-a)
   \boxtimes \Omega^{-b}_{\PP^n}(-b)\otimes F(r,s))=0
\end{equation}
for any $q>0$, $0 \geq a \geq -m$, $0 \geq b \geq -n$, $r+s \geq
-q+1$ and $r,s \geq -q+1$. To this end,  we will first prove the
following claim:

\vskip 2mm \noindent {\bf Claim:} For any $i>0$, $0 \geq b \geq
-n$, $r+s \geq -i$, $r \geq -i$ and $s \geq -i+1$ \[ H^i(
\Omega^{-b}_{\PP^n}(-b)\otimes F(r,s))=0.
\] {\bf Proof of the Claim:} We will prove it by induction on $b$.
By (\ref{HW2}), for $i>0$, $b=-n$,  $r+s \geq -i$, $r \geq -i$ and
$s \geq -i+1$, we have
\[ H^i(\Omega^{n}_{\PP^n}(n)\otimes F(r,s))=H^i(F(r,s-1))=0. \]
Now take $0\le b<-n$ and consider  on $X$ the  exact sequence

\[0 \rightarrow \Omega^{-b+1}_{\PP^n}(-b+1) \otimes F(r,s-1) \rightarrow
\cO_{\PP^n}(-1)^{n+1 \choose -b+1}\otimes F(r,s) \rightarrow
\Omega^{-b}_{\PP^n}(-b)\otimes F(r,s) \rightarrow 0 \]  and the
cohomological exact sequence associated to it
\[ \cdots \rightarrow H^i(F(r,s-1)^{n+1 \choose -b+1}) \rightarrow
H^i(\Omega^{-b}_{\PP^n}(-b)\otimes F(r,s)) \rightarrow
H^{i+1}(\Omega^{-b+1}_{\PP^n}(-b+1) \otimes F(r,s-1))\rightarrow
\cdots .\] Applying (\ref{HW2}) and  hypothesis of induction we
get $H^i(\Omega^{-b}_{\PP^n}(-b)\otimes F(r,s)) =0$ for all $i>0$,
 $r+s \geq -i$, $r \geq -i$ and $s \geq
-i+1$ which finishes the proof of the claim.

\vspace{3mm}

Let us now prove (\ref{2CMR}) by decreasing induction on
$p:=-a-b$, $0 \leq p \leq d$. If $p=d$, then $a=-m$, $b=-n$ and by
(\ref{HW2}) for any $q>0$, $r+s \geq -q+1$ and $r,s \geq -q+1$ we
have
\[ H^q( \Omega^{m}_{\PP^m}(m)
   \boxtimes \Omega^{n}_{\PP^n}(n)\otimes F(r,s))= H^q
   (F(r-1,s-1))=0.\]
   Assume that (\ref{2CMR}) holds for $p+1$ and fix $a$, $b$ such that
   $0<p=-a-b<d$. If $a=0$ or $a=-m$ the result follows from the
   claim.
   So, we can assume $0 > a > -m$ and we consider  on $X$ the exact sequence
\[0 \rightarrow \Omega^{-a+1}_{\PP^m}(-a+1) \boxtimes \Omega^{-b}_{\PP^n}(-b)\otimes F(r-1,s)
 \rightarrow \cO_{\PP^m}(-1)^{m+1 \choose -a+1}\boxtimes \Omega^{-b}_{\PP^n}(-b)\otimes
 F(r,s)\] \[
\rightarrow \Omega^{-a}_{\PP^m}(-a)\boxtimes
\Omega^{-b}_{\PP^n}(-b)\otimes F(r,s) \rightarrow 0 \] and the
cohomological exact sequence associated to it
\[ \cdots   \rightarrow H^q( \cO_{\PP^m}(-1)^{m+1 \choose -a+1}\boxtimes
 \Omega^{-b}_{\PP^n}(-b)\otimes F(r,s))
\rightarrow H^q(\Omega^{-a}_{\PP^m}(-a)\boxtimes
\Omega^{-b}_{\PP^n}(-b)\otimes F(r,s)) \] \[ \rightarrow H^{q+1}(
\Omega^{-a+1}_{\PP^m}(-a+1) \boxtimes
\Omega^{-b}_{\PP^n}(-b)\otimes F(r-1,s)) \rightarrow \cdots. \]

\vskip 2mm \noindent By hypothesis of induction $H^{q+1}(
\Omega^{-a+1}_{\PP^m}(-a+1) \boxtimes
\Omega^{-b}_{\PP^n}(-b)\otimes F(r-1,s))=0$ for any $q>0$, $0 \geq
a \geq -m$, $0 \geq b \geq -n$, $r+s \geq -q+1$ and $r,s \geq
-q+1$ and it follows from the Claim that $H^q(
\cO_{\PP^m}(-1)^{m+1 \choose -a+1}\boxtimes
 \Omega^{-b}_{\PP^n}(-b)\otimes F(r,s))=0$ for any $q>0$, $0 \geq a
\geq -m$, $0 \geq b \geq -n$, $r+s \geq -q+1$ and $r,s \geq -q+1$.
Hence $H^q(\Omega^{-a}_{\PP^m}(-a)\boxtimes
\Omega^{-b}_{\PP^n}(-b)\otimes F(r,s))=0$ for any $q>0$, $0 \geq a
\geq -m$, $0 \geq b \geq -n$, $r+s \geq -q+1$ and $r,s \geq -q+1$
and this finishes the proof of (\ref{2CMR}).

\vskip 2mm  Let us prove the converse. We will prove that
(\ref{1CMR}) implies
\begin{equation} \label{HW3}
H^{i}(X,F(-s,-t))=0 \mbox{ for all } i>0, \quad i \geq s+t-1,\quad
s,t >0.
\end{equation}

First of all, we will prove by induction on $t$ that for all $i
\geq t$
\begin{equation} \label{-1t}
H^{i}(X,F(-1,-t))=0.
\end{equation}

By (\ref{1CMR}), for $t=1$ and $i>0$ we have
\[H^i(X,F(-1,-1))=H^i(X, \Omega^{m}_{\PP^m}(m) \boxtimes \Omega^{n}_{\PP^n}(n)\otimes F)=0.\]
For $t >1$, we consider  on $X$ the exact sequence
\[0 \rightarrow  F(-1,-t)
 \rightarrow  F(-1,-(t-1))^{n+1} \rightarrow \Omega^{n-1}_{\PP^n}(n-t+1)\otimes F(-1,0)
  \rightarrow 0 \] and the
cohomological exact sequence associated to it
\[
\rightarrow H^{i-1}(\Omega^{n-1}_{\PP^n}(n-t+1)\otimes F(-1,0))
\rightarrow H^i(F(-1,-t)) \rightarrow H^i(F(-1,-(t-1))^{n+1})
\rightarrow \cdots . \] By hypothesis of induction
$H^i(F(-1,-(t-1)))=0$ for any $i>t-1$. Hence it is enough to prove
the following Claim:

\vspace{3mm}

 \noindent {\bf Claim 1:} For any $j$, $0 \leq j \leq
t-2$ and $i \geq j+1$,
\[ H^{i}(\Omega^{n-t+1+j}_{\PP^n}(n-t+1)\otimes F(-1,0))=0. \]
{\bf Proof of the Claim 1:} We will proceed by induction on $j$.
For $j=0$ and $i>0$, by (\ref{1CMR})
\[ H^i(\Omega^{n-t+1}_{\PP^n}(n-t+1)\otimes F(-1,0))=
H^i(\Omega^{m}_{\PP^m}(m) \boxtimes
\Omega^{n-t+1}_{\PP^n}(n-t+1)\otimes F)= 0.\] For $0 < j \leq
t-2$, consider the exact sequence on $X$
\[0 \rightarrow \Omega^{n-t+1+j}_{\PP^n}(n-t+1)\otimes F(-1,0)
 \rightarrow  F(-1,-j) ^{{n+1\choose n+1-t+j}}\rightarrow \Omega^{n-t+j}_{\PP^n}(n-t+1)\otimes
 F(-1,0) \rightarrow 0 \] and the
cohomological exact sequence associated to it
\[ \rightarrow H^{i-1} (\Omega^{n-t+j}_{\PP^n}(n-t+1)\otimes
 F(-1,0)) \rightarrow H^{i}(\Omega^{n-t+1+j}_{\PP^n}(n-t+1)\otimes
 F(-1,0)) \] \[ \rightarrow H^i( F(-1,-j)) ^{{n+1\choose n+1-t+j}}\rightarrow \cdots .\]
 By hypothesis of induction $H^{i-1} (\Omega^{n-t+j}_{\PP^n}(n-t+1)\otimes
 F(-1,0))=0$ for any $i>j$ and since $j<t$, by hypothesis of induction on $t$,
 $H^i(F(-1,-j))=0$ for any $i \geq j$. Thus $H^{i}(\Omega^{n-t+1+j}_{\PP^n}(n-t+1)\otimes
 F(-1,0))=0$ for any $i \geq j+1$. This finishes the proof of
 the Claim 1 and the proof of (\ref{-1t}).

  By symmetry, for all $i \geq t$ we also have
\begin{equation} \label{sim}
H^{i}(X,F(-t,-1))=0.
\end{equation}

 \vspace{2mm}

Now, by induction on $s>0$ we will prove that (\ref{HW3}) holds
for any $t>0$, $i \geq t+s-1$. The case $s=1$ is already done.
Take $s>1$  and we will see that  \begin{equation} \label{fi}
H^i(F(-s,-t))=0 \mbox{  for any } t\le i, \quad  i \geq s+t-1.
\end{equation} To this end, we will prove that the following cohomology
groups vanish on $X$:
\begin{itemize}
\item[(a)] For any $i \geq s \geq 1$,
\[H^i( \Omega^{n-t+1}_{\PP^n}(n-t+1)\otimes F(-s,0))=0. \]
\item[(b)] For any $j$, $1 \leq j \leq t-1$, $t>0$ and $i \geq s+t-1$,
\[H^i(\Omega^{n-t+j}_{\PP^n}(n-t+1)\otimes F(-s,0))=0. \]
\end{itemize}

\noindent $(a)$ We proceed by induction on $s$. For $s=1$, by
(\ref{1CMR}), for any $i>0$
\[H^i( \Omega^{n-t+1}_{\PP^n}(n-t+1)\otimes F(-1,0))=
H^i(\Omega^{m}_{\PP^m}(m) \boxtimes
\Omega^{n-t+1}_{\PP^n}(n-t+1)\otimes F)=0. \]

Fix $s>1$ and let us prove

\vspace{3mm}

 \noindent  {\bf Claim 2:} For any $j$, $1 \leq
j \leq s-1$ and $i \geq j$,
\[H^i(\Omega^{m-s+j}_{\PP^m}(m-s+1) \boxtimes \Omega^{n-t+1}_{\PP^n}(n-t+1)\otimes F)=0. \]

\noindent {\bf Proof of Claim 2:} For $j=1$, by (\ref{1CMR}) \[
H^i(\Omega^{m-s+1}_{\PP^m}(m-s+1) \boxtimes
\Omega^{n-t+1}_{\PP^n}(n-t+1)\otimes F)=0. \] For $1 <j \leq s-1$
consider the exact sequence on $X$
\[0 \rightarrow \Omega^{m-s+j}_{\PP^m}(m-s+1) \boxtimes \Omega^{n-t+1}_{\PP^n}(n-t+1)
 \otimes F  \rightarrow \Omega^{n-t+1}_{\PP^n}(n-t+1)\otimes
 F(-j+1,0)^{{m+1\choose m-s+j}}\] \[
\rightarrow \Omega^{m-s+j-1}_{\PP^m}(m-s+1) \boxtimes
\Omega^{n-t+1}_{\PP^n}(n-t+1)\otimes F \rightarrow 0
\] and the cohomological exact sequence associated to it
 \[\cdots \rightarrow H^{i-1}(\Omega^{m-s+j-1}_{\PP^m}(m-s+1) \boxtimes
\Omega^{n-t+1}_{\PP^n}(n-t+1)\otimes F)  \rightarrow \] \[
H^i(\Omega^{m-s+j}_{\PP^m}(m-s+1) \boxtimes
\Omega^{n-t+1}_{\PP^n}(n-t+1)
 \otimes F) \rightarrow H^i(\Omega^{n-t+1}_{\PP^n}(n-t+1)\otimes
 F(-j+1,0))^{{m+1\choose m-s+j}}\rightarrow \cdots . \]
By hypothesis of induction on $j$,
$H^{i-1}(\Omega^{m-s+j-1}_{\PP^m}(m-s+1) \boxtimes
\Omega^{n-t+1}_{\PP^n}(n-t+1)\otimes F)=0$ and by hypothesis of
induction on $s$, for $j-1 < s$,
$H^i(\Omega^{n-t+1}_{\PP^n}(n-t+1)\otimes
 F(-j+1,0))=0$. Thus $H^i(\Omega^{m-s+j}_{\PP^m}(m-s+1) \boxtimes
\Omega^{n-t+1}_{\PP^n}(n-t+1) \otimes F)=0$ and this finishes the
proof of Claim 2.

Using the cohomological exact sequence on $X$
\[\rightarrow H^{i-1}(\Omega^{m-1}_{\PP^m}(m-s+1) \boxtimes
\Omega^{n-t+1}_{\PP^n}(n-t+1)\otimes F)  \rightarrow H^i(
\Omega^{n-t+1}_{\PP^n}(n-t+1)
 \otimes F(-s,0))  \] \[ \rightarrow H^i(\Omega^{n-t+1}_{\PP^n}(n-t+1)\otimes
 F(-s+1,0))^{m+1}\rightarrow \cdots  \]
  by hypothesis of induction and by Claim 2, we get that for any $i \geq s$, $s \geq 1$, $H^i(
\Omega^{n-t+1}_{\PP^n}(n-t+1)\otimes F(-s,0))=0$ which finishes
the proof of $(a)$.

\vspace{3mm}

$(b)$ The case $j=1$ follows from $(a)$. Fix $j>0$ and consider
the cohomological exact sequence
\[\rightarrow H^{i-1}(\Omega^{n-t+j-1}_{\PP^n}(n-t+1)\otimes F(-s,0))  \rightarrow H^i(
\Omega^{n-t+j}_{\PP^n}(n-t+1)
 \otimes F(-s,0))  \] \[ \rightarrow H^i( F(-s,-j+1))^{{n+1\choose n-t+j}}\rightarrow \cdots  \]
associated to the exact sequence on $X$
\[ 0 \rightarrow \Omega^{n-t+j}_{\PP^n}(n-t+1)
 \otimes F(-s,0) \rightarrow  F(-s,-j+1)^{{n+1\choose n-t+j}} \rightarrow
  \Omega^{n-t+j-1}_{\PP^n}(n-t+1)\otimes F(-s,0) \rightarrow 0.\]
Since $j-1 < t$, by the first cases $H^i( F(-s,-j+1))=0$ and by
hypothesis of induction
$H^{i-1}(\Omega^{n-t+j-1}_{\PP^n}(n-t+1)\otimes F(-s,0))=0$. Thus,
for any $j$, $1 \leq j \leq t-1$ and $i \geq s+t-1$,
$H^i(\Omega^{n-t+j}_{\PP^n}(n-t+1)\otimes F(-s,0))=0$ which
finishes the proof of $(b)$.

\vspace{3mm} Finally, to prove (\ref{fi}), we proceed by induction
on $t>0$. By (\ref{sim}), the case $t=1$ is already done, so we
fix $t>1$ and we consider the following exact sequence on $X$
\[0 \rightarrow F(-s,-t) \rightarrow F(-s,-t+1)^{n+1} \rightarrow
\Omega^{n-1}_{\PP^n}(n-t+1)\otimes F(-s,0) \rightarrow 0 \] and
the associated cohomological exact sequence
\[ \rightarrow H^{i-1}(\Omega^{n-1}_{\PP^n}(n-t+1)\otimes F(-s,0)) \rightarrow
H^i(F(-s,-t)) \rightarrow H^i(F(-s,-t+1))^{n+1} \rightarrow \cdots
.
\] By hypothesis of induction $H^i(F(-s,-t+1))=0$ and it follows
from $(b)$ that $H^{i-1}(\Omega^{n-1}_{\PP^n}(n-t+1)\otimes
F(-s,0))=0$. Hence, $H^i(F(-s,-t))=0$ for any  $t>0$, $i \geq
s+t-1$ and this proves what we want.
\end{proof}


\vskip 4mm As a consequence of this last Theorem we have,

\begin{corollary} Let $F$ be a coherent sheaf on $X=\PP^m\times
\PP^n$. Set $d=n+m$.

(1) If $F$   is $p$-regular with respect to $\cB$ in the sense of
Definition \ref{new} and $p=\lambda(d+1)+\rho$, $0<\rho\le d+1$,
then $F$  is $((\lambda +2)(m+1),(\lambda +2)(n+1))$-regular in
the sense of Hoffman and Wang; and

(2) If $F$  is $(s,r)$-regular in the sense of Hoffman and Wang
with $r=\lambda(m+1)+t$, $0<t\le m+1$, and  $s=\mu(m+1)+x$,
$0<x\le n+1$ then $F$ is $(max(\lambda ,\mu )(d+1)+1)$-regular
with respect to $\cB$ in the sense of Definition \ref{new}.
\end{corollary}
\begin{proof} (1) By Proposition \ref{sup}, if $F$   is
$p=\lambda(d+1)+\rho$-regular with respect to $\cB$, then $F$ is
$((\lambda +2)(d+1)-d)$-regular with respect to $\cB $ as well. By
Lemma \ref{key1}, $F((\lambda +2)(m+1),(\lambda +2)(n+1))$ is
$(-d)$-regular with respect to $\cB $. Applying Theorem
\ref{mainthm2}, we get that $F((\lambda +2)(m+1),(\lambda
+2)(n+1))$   is $(0,0)$-regular in the sense of Hoffman and Wang
and hence $F$   is $((\lambda +2)(m+1),(\lambda +2)(n+1))$-regular
in the sense of Hoffman and Wang.

\vskip 2mm (2) Set $\phi =max(\lambda,\mu)$. By \cite{HW};
Proposition 2.7, $F$  is $((\phi +1)(m+1),(\phi +1)(n+1))$-regular
in the sense of Hoffman and Wang. Therefore, $F((\phi
+1)(m+1),(\phi +1)(n+1))$ is (0,0)-regular in the sense of Hoffman
and Wang and applying Theorem \ref{mainthm2} we obtain that
$F((\phi +1)(m+1),(\phi +1)(n+1))$ is $(-d)$-regular with respect
to $\cB $ and so $F$ is $(\phi(d+1)+1)$-regular with respect to
$\cB$ in the sense of Definition \ref{new}.
\end{proof}


\section{Final remark and open problem}

The notion of regularity that we have introduced in section 4
applies to any coherent sheaf on a large class of smooth
projective varieties: projective spaces, multiprojective spaces,
hyperquadric varieties, Grassmannians, etc. More precisely, it
applies to coherent sheaves on any $n$-dimensional smooth
projective variety which has an $n$-block collection $\cB=( \cE_0,
\cE_1, \cdots ,\cE_n)$ of type $(\alpha _0, \alpha _1, \cdots
,\alpha _n)$ of coherent sheaves on $X$ which generates the
derived category of bounded complexes $\cD = D^b(\cO_X-mod)$.
Hence, we are led to pose the following question/problem:

\begin{problem} To characterize $n$-dimensional smooth projective varieties
 which have an
$n$-block collection $\cB=( \cE_0, \cE_1, \cdots ,\cE_n)$,
$\cE_j=(E_1^j, E_2^j,\cdots, E_{\alpha _j}^j)$ of coherent sheaves
on $X$ which generates $\cD$.
\end{problem}

\vskip 4mm

\noindent {\bf Note added in Proof:}

\noindent  Problem \ref{prob1}  is closely related to Dubrovin's
conjecture concerning the semisimplicity of the quantum cohomology
algebra. More precisely, it states

\vspace{3mm}

\noindent {\bf Conjecture:} \rm (Dubrovin \cite{Du}; Conjecture
4.2.2 (1)) Let $X$ be a smooth complex compact variety.
 The even quantum
cohomology ring of $X$ is generically semisimple if and only if
$X$ is a Fano variety and the category $\cD$ admits a full
exceptional collection of length equal to $\sum_q H^{q,q}(X)$.


\end{document}